\documentclass[preprint,1p,leqno]{elsarticle}
\usepackage{amsmath}
\usepackage{amssymb}
\usepackage{mathrsfs}
\usepackage{bm}
\usepackage{arydshln}

\usepackage{color}
\usepackage[svgnames]{xcolor}
\usepackage{marvosym}
\usepackage{ulem}

\usepackage[pagewise]{lineno}
\modulolinenumbers[1]

\usepackage[]{graphicx}
\usepackage{psfrag}
\usepackage{pdfpages}
\usepackage[]{xcolor}

\usepackage[makeroom]{cancel}

\usepackage{amsthm}

\usepackage{bigstrut}
\usepackage{multirow}
\usepackage{subfig}
\usepackage{easyReview}
\setreviewson
\usepackage{hyperref}

\usepackage{etoolbox}
\apptocmd{\thebibliography}{\setlength{\itemsep}{0pt}}{}{}

\usepackage{geometry}
\newcommand{\bu}{\boldsymbol{u}}
\newcommand{\bF}{\boldsymbol{f}}
\newcommand{\bw}{\boldsymbol{\omega}}

\newcommand{\bn}{\boldsymbol{n}}
\newcommand{\bj}{\boldsymbol{j}}

\newcommand{\bB}{\boldsymbol{B}}
\newcommand{\bH}{\boldsymbol{H}}
\newcommand{\bE}{\boldsymbol{E}}

\newcommand{\bv}{\boldsymbol{v}}

\newcommand{\Rn}{R_f}
\newcommand{\Rm}{R_m}

\journal{Somewhere}
\begin{document}
	\begin{frontmatter}
		
		\title{A dual-field structure-preserving mixed finite element discretization for incompressible Hall MHD equations}
		
		\author[1,2]{Yi Zhang\corref{cor1}}
		\ead{zhangyi_aero@hotmail.com}

		\affiliation[1]{organization={School of Mathematics and Computing Science, Guilin University of Electronic Technology},
			city={Guilin},
			country={China}}
		
		\affiliation[2]{organization={Center for Applied Mathematics of Guangxi (GUET)},
			city={Guilin},
			country={China}}
		
		\cortext[cor1]{Corresponding author}

		%\fntext[fn1]{This is the first author footnote.}
		
		\begin{abstract}
			In this paper, a novel dual-field structure-preserving mixed finite element discretization for the incompressible Hall magnetohydrodynamics (MHD) equations is introduced. It preserves pointwise mass conservation, magnetic Gauss law, and conservation of current density. It also satisfies a discrete energy law that exactly captures the correct energy dissipation rate in the dissipative regime and recovers energy conservation in the ideal limit. Numerical experiments are presented to verify the temporal and spatial accuracy and to illustrate the structure-preserving properties of the scheme.
		\end{abstract}
	\end{frontmatter}
	
	%_____________________________________________________________________________________________________
	%>>>>>>>>>>>>>>>>>>>>>>>>>>>>>>>>>>>>>> SECTION <<<<<<<<<<<<<<<<<<<<<<<<<<<<<<<<<<<<<<<<<<<<<<<<<<<<<<<
	\section{Introduction} \label{Sec: Introduction}
	
	\subsection{The Hall MHD model} \label{sub: Hall MHD}
	In $\mathbb{R}^3$, we consider a bounded, contractible subdomain $\Omega$ with a Lipschitz boundary $\partial\Omega$. 
	In a space-time domain, $\Omega\times\left(0,T \right] $, a dimensionless form of three-dimensional incompressible Hall magnetohydrodynamics (MHD) equations, hereinafter referred to as the Hall MHD model, is written as,  \cite{LAAKMANN2023112410, doi:10.1137/23M1553844}, 
	\begin{subequations}\label{eq: mixed Hall MHD}
		\begin{align}
			\partial _{t}\bu +\bw\times \bu  + \Rn^{-1}\nabla\times \bw  - \mathsf{c}\bj \times\bB  + \nabla  P  &= \bF\quad &&\text{in }\Omega\times\left(0,T \right], \label{eq: mixed Hall MHD a}\\
			\bw &=\nabla\times\bu\quad &&\text{in }\Omega\times\left(0,T \right],\label{eq: mixed Hall MHD b}\\
			\nabla\cdot   \bu  &=0\quad &&\text{in }\Omega\times\left(0,T \right],\label{eq: mixed Hall MHD c}\\
			\bj &=\nabla\times\boldsymbol{B} \quad &&\text{in }\Omega\times\left(0,T \right],\label{eq: mixed Hall MHD d}\\
			\partial_{t}\bB + \nabla\times \bE  &= \boldsymbol{0} \quad &&\text{in }\Omega\times\left(0,T \right], \label{eq: mixed Hall MHD e}\\
			\Rm^{-1}\bj - \left( \bE +\bu\times\bB\right) + \mathsf{h}\bj\times\bB &= \boldsymbol{0}\quad &&\text{in }\Omega\times\left(0,T \right],\label{eq: mixed Hall MHD f} 
		\end{align}
	\end{subequations}
	where $\partial _{t}:=\frac{\partial}{\partial t}$. It governs the dynamics of fluid velocity $\bu$, vorticity $\bw$, electric current density $\bj$, magnetic flux density $\bB$, total pressure $P:=\frac{1}{2}\bu\cdot\bu + p$ (with $p$ being the static pressure), and electric field strength $\bE$, subject to an external body force $\bF$, initial conditions $\bu(\boldsymbol{x},0) = \bu^{0}$ and $\bB(\boldsymbol{x},0) =\bB^{0}$, and proper boundary conditions. Four dimensionless parameters are the fluid Reynolds number $\Rn$, the coupling number for Lorentz force $\mathsf{c}$, the magnetic Reynolds number $\Rm$, and the strength factor of Hall effect $\mathsf{h}$. It is clear that when $\mathsf{h}=0$, \eqref{eq: mixed Hall MHD} reduces to incompressible MHD equations. 
	
	The Hall MHD model arises in a wide range of physical phenomena such as magnetic reconnection, space plasmas, and fusion devices \cite{Birn2001GEM,BirnPriest2007}. This highly nonlinear model involves two subsystems, namely the incompressible Navier-Stokes equations, \eqref{eq: mixed Hall MHD a} - \eqref{eq: mixed Hall MHD c}, and Maxwell equations with Hall effect, \eqref{eq: mixed Hall MHD d} - \eqref{eq: mixed Hall MHD f}, which are nonlinear because of the fluid convection term, $\bw\times\bu$ in \eqref{eq: mixed Hall MHD a}, and the Hall effect term, $\mathsf{h}\bj\times\bB$ in \eqref{eq: mixed Hall MHD f}. They are coupled through additional nonlinear terms, $\mathsf{c}\bj\times \bB $ in \eqref{eq: mixed Hall MHD a} and $\bu\times\bB$ in \eqref{eq: mixed Hall MHD f}, which represent Lorentz force and electromotive force, respectively. Compared with the MHD model, the extra nonlinear term, i.e. the Hall effect term, further strengthens the interconnection between variables and makes the flow fully three-dimensional. This can be seen from \eqref{eq: mixed Hall MHD f}; if the flow is two-dimensional, the Hall effect term will contribute as a vector field while other terms represent scalar fields, which is contradictory. These complexities render the Hall MHD model a challenging problem to analyze and simulate numerically.
	
	The Maxwell part in the Hall MHD model can also be written as the magnetic induction (or Helmholtz) equation with the Hall effect term, i.e.,
	\begin{equation}\label{eq: eB}
		\partial_{t}\bH + \Rm^{-1}\Delta\bH - \nabla\times\left( \bu\times\bH\right) +\mathsf{h}\nabla\times\left(\left(\nabla\times\bH \right) \times\bH \right) =\boldsymbol{0},
	\end{equation}
	where $\bH$ is the magnetic field strength, which, in a dimensionless form, satisfies\footnote{In differential forms, they are a pair of dual forms which are connected by the Hodge-star operator. See \ref{App: df HMHD}.}
	\begin{equation}\label{eq: BH}
		\bB=\bH.
	\end{equation}
	If we combine \eqref{eq: mixed Hall MHD d} - \eqref{eq: mixed Hall MHD f} and eliminate $\bE$ and $\bj$, we can get the same magnetic induction equation \eqref{eq: eB} but for $\bB$, which also implies \eqref{eq: BH}. Also see (2) of \cite{doi:10.1137/23M1553844}. 
	
	In the present work, we group \eqref{eq: mixed Hall MHD} and \eqref{eq: eB} and work on a system written as 
	\begin{subequations}\label{eq: mixed Hall MHD BH}
		\begin{align}
			\partial _{t}\bu +\bw\times \bu  + \Rn^{-1}\nabla\times \bw  - \mathsf{c}\bj \times\bH  + \nabla  P  &= \bF\quad &&\text{in }\Omega\times\left(0,T \right], \label{eq: mixed Hall MHD HB a}\\
			\bw &=\nabla\times\bu\quad &&\text{in }\Omega\times\left(0,T \right],\label{eq: mixed Hall MHD HB b}\\
			\nabla\cdot   \bu  &=0\quad &&\text{in }\Omega\times\left(0,T \right],\label{eq: mixed Hall MHD HB c}\\
			\bj &=\nabla\times\boldsymbol{B} \quad &&\text{in }\Omega\times\left(0,T \right],\label{eq: mixed Hall MHD HB d}\\
			\partial_{t}\bB + \nabla\times \bE  &= \boldsymbol{0} \quad &&\text{in }\Omega\times\left(0,T \right], \label{eq: mixed Hall MHD HB e}\\
			\Rm^{-1}\bj - \left( \bE +\bu\times\bH\right) + \mathsf{h}\bj\times\bH &= \boldsymbol{0}\quad &&\text{in }\Omega\times\left(0,T \right],\label{eq: mixed Hall MHD HB f} \\
			\partial_{t}\bH + \Rm^{-1}\Delta\bH - \nabla\times\left( \bu\times\bH\right) +\mathsf{h}\nabla\times\left(\left(\nabla\times\bH \right) \times\bB \right) &=\boldsymbol{0}\quad &&\text{in }\Omega\times\left(0,T \right],\label{eq: mixed Hall MHD HB g}
		\end{align}
	\end{subequations}
	where we have used \eqref{eq: BH} to replace $\bB$ by $\bH$ in nonlinear terms $\mathsf{c}\bj \times\bH $, $\bu\times\bH $, $\mathsf{h}\bj\times\bH$, and to replace the second $\bH$ in term $ \mathsf{h}\nabla\times\left(\left(\nabla\times\bH \right) \times\bB \right) $ by $\bB$. % This mixture is one of the keys to linearizing the present discretization, as we will explain later in Section~\ref{Sec: formulation} and Section~\ref{Sec: discretization}. 
	
	It seems that we have complicated the formulation by adding one additional variable and one equation, namely $\bH$ and \eqref{eq: mixed Hall MHD HB g}, compared to \eqref{eq: mixed Hall MHD}. But we will show that, under a carefully designed temporal discretization to a weak form of \eqref{eq: mixed Hall MHD BH}, we can linearize the system while preserving some important structures of it. The formulation \eqref{eq: mixed Hall MHD BH} is mainly inspired by the dual-field method initially proposed for incompressible Navier-Stokes equations in \cite{ZHANG2022110868} and then extended to the MHD model in \cite{MAO2025114130}, where dual representations of solutions are sought. Taking the velocity field $\bu$ in \cite{ZHANG2022110868} as an example, under a mixed finite element framework, instead of choosing between seeking $\bu_h\in H(\mathrm{curl})$ or seeking $\bu_h\in H(\mathrm{div})$, the authors decide to solve for two variables, $\bu^1_h\in H(\mathrm{curl})$ and $\bu^2_h\in H(\mathrm{div})$, i.e., the dual representations of $\bu_h$, in dual systems. By doing so, the nature of $\bu\in H(\mathrm{curl}) \cap H(\mathrm{div})$ is captured in an implicit sense. And the dual systems can be linearized by borrowing information from each other and are then solved in a leapfrog-type temporal scheme. With the present formulation \eqref{eq: mixed Hall MHD BH}, having evolution equations for both $\bB$ and $\bH$, we are able to build a dual-field discretization for the magnetic field, which in turn leads to a linearized structure-preserving scheme for the incompressible Hall MHD model. 
	The relation between $\bB$ and $\bH$ and other structures behind \eqref{eq: mixed Hall MHD BH} is clearer in differential forms. See \ref{App: df HMHD} for a brief explanation. 
	
	Compared to existing dual-field methods \cite{ZHANG2022110868, MAO2025114130}, the present method has an obvious difference that it uses a strong form having dual variables present already, i.e. $\bB$ and $\bH$ in \eqref{eq: mixed Hall MHD BH}. While the existing dual-field methods use strong forms without dual variables, see \cite[(9)]{ZHANG2022110868} and \cite[(2.2)]{MAO2025114130}, and the sense of being dual-field comes later with the spatial discretization employing two mixed weak forms that are dual to each other. We will come back to this point when the weak form of the present work is introduced in Section~\ref{Sec: formulation}.
	
	In this work, we will restrict the mixed weak form and its analysis to homogeneous boundary conditions,
	\begin{equation}\label{eq: bc}
		P = 0,\quad \boldsymbol{u}\times\boldsymbol{n} = \boldsymbol{0},\quad \bB\times\bn = \bH\times\bn = \boldsymbol{0} \quad \text{on }\partial\Omega\times(0,T].
	\end{equation}
	And it will also be shown that it is possible to extend the mixed weak form to general boundary conditions. Finally, to close the problem, intial conditions, $\bu^{0}$, $\bB^{0}$ and $\bH^{0}(=\bB^{0})$, are supplemented to \eqref{eq: mixed Hall MHD BH}. 
	
	In the incompressible Hall MHD equations, the divergence-free condition of $\bu$, i.e., \eqref{eq: mixed Hall MHD HB c}, implies the incompressibility of the fluid. And if the initial condition $\bB^{0}$ satisfies $\nabla\cdot\bB^0 = 0$, the evolution equation of $\bB$, \eqref{eq: mixed Hall MHD HB e}, will ensure that $\nabla\cdot\bB = 0$ in $\Omega\times(0,T]$, which is magnetic Gauss's law. From \eqref{eq: mixed Hall MHD HB d} and the property $\nabla\cdot\nabla\times(\cdot) \equiv 0 $, we know that $\nabla\cdot\bj = 0 $, which represents conservation of current density derived from the law of charge conservation. Furthermore, under boundary conditions in \eqref{eq: bc}, the Hall MHD model satisfies an energy law,
	\begin{equation}\label{eq: energy law}
		\partial_{t}\mathcal{E} = \partial_{t}\left( \mathcal{K}+\mathcal{M}\right)  =\underbrace{-\Rn^{-1}\int_{\Omega}\bw\cdot\bw\ \mathrm{d}\Omega }_{\text{(i)}}\  \underbrace{-\mathsf{c}\Rm^{-1}\int_{\Omega}\bj\cdot\bj\ \mathrm{d}\Omega}_{\text{(ii)}} + \underbrace{\int_{\Omega}\bF\cdot\bu\ \mathrm{d}\Omega}_{\text{(iii)}},
	\end{equation}
	where the terms at the right-hand side represent (i) the dissipation rate due to fluid viscosity, (ii) the dissipation rate due to magnetic resistance, and (iii) the rate of contribution of the external body force. And 
	\[
	\mathcal{K}:=\frac{1}{2}\int_{\Omega}\bu\cdot\bu\ \mathrm{d}\Omega,\quad \mathcal{M}: = \frac{\mathsf{c}}{2}\int_{\Omega}\bB\cdot\bB\ \mathrm{d}\Omega
	\]
	are the (total) fluid kinetic energy and the (total) magnetic energy, respectively. And $\mathcal{E}:=\mathcal{K}+\mathcal{M}$ is called the total energy. 
	In the ideal limit, i.e. $\Rn,\ \Rm \to +\infty$ and $\boldsymbol{f} = \boldsymbol{0}$, this energy law obviously implies conservation of energy, 
	\[
	\partial_{t}\mathcal{E} = 0.
	\]
	
	\subsection{A brief literature survey on structure-preserving methods towards the Hall MHD model}
	The MAC method, originated from the work of Harlow and Welch \cite{harlow1965numerical}, is usually referred to as the first structure-preserving method. Then the mixed finite element method \cite{raviart606mixed,nedelec,boffi2013mixed,arnold2014periodic} opened up a new era of structure-preserving methods in the family of finite element methods. More recently, structure-preserving frameworks such as the mimetic framework introduced by Bochev and Hyman \cite{bochev2006principles} and the discrete exterior calculus (DEC) \cite{hirani2003discrete} have also received successes. They use differential forms instead of the classic vector calculus as the mathematical language, respect geometric representations behind physical variables, and have a tight link to Tonti's discovery that reveals the analogies between differential geometry and algebraic topology in the modeling of physics \cite{tonti2013mathematical}. Having a similar flavour, the finite element exterior calculus (FEEC) \cite{arnold2006finite,arnold2010finite}, bridging the mixed finite element method to differential forms, has also become a powerful tool for structure-preserving discretizations and has furthermore promoted the popularity of the mixed finite element method in the field of structure-preserving methods.  
	
	The benefit on numerical stability of structure-preserving numerical methods for MHD problems was first exposed in the methods conserving magnetic Gauss's law \cite{BRACKBILL1980426}. Afterwards, for a long time, satisfying magnetic Gauss's law at the discrete level exactly becomes a major track in this field, and various methods have been developed which can be classified into, for example,  potential-based methods \cite{hiptmair2018fully,DING2024}, divergence-cleaning methods \cite{TOTH2000605,Balsara_2004}, constrained transport methods \cite{Evans1988SimulationOM,refId0}, Lagrange multiplier (or augmented) methods \cite{Heister2017,10.1093/imanum/drad005}.
	
	Besides magnetic Gauss's law, other structures embedded in the MHD model, like conservation of mass (pointwise incompressibility of the fluid), conservation of current density, the energy law (conservation of energy in the ideal limit or unconditional energy stability in a weaker sense), conservation of (cross- and magnetic-) helicity, are also of great interest. And the mixed finite element method is undoubtedly a very successful method for designing numerical schemes that preserve multiple structures of the MHD model simultaneously. For example, Hu et al. \cite{HU2021110284} have proposed a mixed finite element method that preserves cross- and magnetic-helicity, energy, and magnetic Gauss's law for incompressible MHD equations. For incompressible MHD of a variable fluid density, Gawlik and Gay-Balmaz have introduced a mixed finite element method that preserves energy, cross-helicity (when the fluid density is constant), and magnetic-helicity, total squared density, pointwise incompressibility, and magnetic Gauss's law in \cite{GAWLIK2022110847}. These methods usually lead to big nonlinear systems which are expensive to solve. To linearize the problem, Zhang and Su \cite{ZHANG202345} introduce a structure-preserving discretization which decouples the problem using the mixed finite element method and the FEEC for the incompressible MHD model. This method has a first-order temporal accuracy. A temporally second-order dual-field mixed finite element method that linearizes the problem and satisfies conservation of mass, magnetic Gauss's law, conservation of current density, conservation of energy, and conservation of magnetic and fluid helicity is proposed in \cite{MAO2025114130}. 
	
	A technique that decouples and linearizes the problem while keeping a second-order temporal accuracy is the scalar auxiliary variable (SAV) method \cite{shen2018scalar,shen2019new}. Implemented together with the mixed finite element method, the SAV method has received sucesses for incompressible Navier-Stokes equations \cite{lin2019numerical} and Navier-Stokes equations related coupled problems like the Navier-Stokes-Poisson-Nernst-Planck (NSPNP) system \cite{pan2024linear} and the MHD problem \cite{ZHANG2022110752}. These implementations usually possess an unconditional (modified) energy stability instead of a discrete energy law (or conservation of energy in the ideal limit) because of the introduction of the auxiliary variable. 
	
	In part due to the numerical challenges posed by the Hall MHD model, as discussed in the previous subsection, the development of structure-preserving discretizations for the Hall MHD model has emerged only recently. Representative contributions are the work of Laakmann et al. \cite{LAAKMANN2023112410} and the work of Guo et al. \cite{doi:10.1137/23M1553844}. They both employ the finite element method, and the latter also uses the SAV method. Their differences from the method presented in this paper will be clarified in the coming subsection.
	
	\subsection{Objective and outline of the paper}
	This paper reports a new linear dual-field mixed finite element discretization of the Hall MHD model. The method is structure-preserving; it satisfies pointwise conservation of mass, magnetic Gauss's law, and conservation of current density. Additionally, it preserves a canonical discrete energy law by construction. This method differs from the existing structure-preserving methods for the Hall MHD model \cite{LAAKMANN2023112410,doi:10.1137/23M1553844} in several key aspects, offering both advantages and disadvantages, as follows.
	\begin{itemize}
		\item  The present method is simpler in the sense that it employs a dual-field formulation and only uses the mixed finite element method under a carefully designed temporal scheme, in addition to no other techniques.
		\item  Compared to the method in \cite{LAAKMANN2023112410}, the present method is linear and does not depend on projections. And it by construction preserves an exact discrete energy law instead of an unconditional energy stability as in \cite{doi:10.1137/23M1553844}.
		\item However, the present method does not conserve a hybrid helicity as \cite{LAAKMANN2023112410} does, see (3.7), (3.8) therein. And it does not decouple the problem as \cite{doi:10.1137/23M1553844} does, such that the computation will be more expensive.
	\end{itemize}
	
	The remainder of the paper is organized as follows. In Section~\ref{Sec: formulation}, we introduce preliminaries, the semi-discrete formulation and its properties. Its fully discrete version and the properties follow in Section~\ref{Sec: discretization}. Supportive numerical experiments are presented in Section~\ref{Sec: numerical}. Finally, Section~\ref{Sec: conclusions} summarizes the paper and discusses limitations and future directions of the present work.
	
	\section{A spatially discrete weak form} \label{Sec: formulation}
	\subsection{Function spaces} \label{sub: spaces}
	Let $\left\langle\cdot,\cdot\right\rangle_{\Omega}$ denote the $L^2$-inner product over the domain $\Omega$. The well-known Hilbert complex of Sobolev spaces is 
	\[
	\mathbb{R}\hookrightarrow H^1(\Omega)\stackrel{\nabla}{\longrightarrow}\boldsymbol{H}(\Omega;\mathrm{curl})\stackrel{\nabla\times}{\longrightarrow}\boldsymbol{H}(\Omega;\mathrm{div})\stackrel{\nabla\cdot}{\longrightarrow}L^2(\Omega)\longrightarrow 0 .
	\]
	The finite-dimensional spaces used in this work are
	\[
	G\subset H^1(\Omega),\ \boldsymbol{C}\subset \boldsymbol{H}(\Omega;\mathrm{curl}),\ \boldsymbol{D}\subset \boldsymbol{H}(\Omega;\mathrm{div}),\ S\subset L^2(\Omega),
	\]
	which also form a Hilbert complex,
	\begin{equation}\label{Eq: discrete complex}
		\mathbb{R}\hookrightarrow G\stackrel{\nabla}{\longrightarrow}\boldsymbol{C}\stackrel{\nabla\times}{\longrightarrow}\boldsymbol{D}\stackrel{\nabla\cdot}{\longrightarrow} S \longrightarrow 0 .
	\end{equation}
	And we will use $\boldsymbol{C}_{0}$ to denote the subspace $\boldsymbol{C}_{0}:=\left\lbrace\boldsymbol{a}_{h}\left|\boldsymbol{a}_{h}\in\boldsymbol{C}, \boldsymbol{a}_{h}\times\boldsymbol{n}=\boldsymbol{0} \text{ on }\partial\Omega\right.\right\rbrace$ where $\boldsymbol{n}$ is the outward unit normal vector.
	To adapt to the nonlinear terms of the Hall MHD model, the finite-dimensional spaces should also satisfy one additional regularity that the $L^2$-integrability of the trilinear term,
	\begin{equation}\label{Eq:regularity}
		\mathcal{A}(\boldsymbol{\alpha}_{h}, \boldsymbol{\beta}_{h}, \boldsymbol{\gamma}_{h}) := \left\langle\boldsymbol{\alpha}_{h}\times\boldsymbol{\beta}_{h},\boldsymbol{\gamma}_{h}\right\rangle_{\Omega},\ \forall \boldsymbol{\alpha}_{h}, \boldsymbol{\beta}_{h}, \boldsymbol{\gamma}_{h} \in \boldsymbol{C} \text{ or } \boldsymbol{D},
	\end{equation}
	is guaranteed. Note that $\mathcal{A}$ is skew-symmetric with respect to any two of the three entries.
	
	\subsection{The weak form}
	In this subsection, we present a spatially discrete weak form of the Hall MHD model. Its propesties are analyzed in coming subsections. For a concise presentation, hereinafter, we call a spatially discrete object a semi-discrete one and omit the temporal part, i.e., $\times(0, T]$, in the notations of space-time domains. 
	
	We consider the following semi-discrete weak form of the mixed system \eqref{eq: mixed Hall MHD BH} under homogeneous boundary conditions \eqref{eq: bc}. Given $\bF\in \left[  L^2(\Omega)\right] ^3$, seek $\left(\bu_{h},\bw_{h}, P_{h}, \bE_{h},\bB_{h}, \bj_{h},\bH_{h}\right)\in \boldsymbol{D}\times \boldsymbol{C}\times S\times \boldsymbol{C}\times \boldsymbol{D}\times \boldsymbol{C}\times\boldsymbol{C}_{0}$, such that 
	\begin{subequations}\label{Eq: total system}
		\begin{align}
			\left\langle \partial _{t}\bu_{h},\bv_{h}\right\rangle _{\Omega} +\mathcal{A}\left( \bw_{h}, \bu_{h},\bv_{h}\right)    + \Rn^{-1}\left\langle \nabla\times \bw_{h} ,\bv_{h}\right\rangle _{\Omega} \qquad \label{eq: a}
			\\
			- \mathsf{c}\mathcal{A}\left(  \bj_{h} ,{\bH_{h}},\bv_{h}\right)  - \left\langle P_{h},\nabla\cdot\bv_{h}\right\rangle _{\Omega}  &= \left\langle \bF,\bv_{h}\right\rangle _{\Omega}&&\forall \bv_{h}\in \boldsymbol{D},\nonumber\\
			\left\langle \bw_{h},\boldsymbol{w}_{h} \right\rangle _{\Omega}- \left\langle\bu_{h}, \nabla\times\boldsymbol{w}_{h}\right\rangle _{\Omega} &=0&&\forall \boldsymbol{w}_{h}\in \boldsymbol{C},\label{eq: b}
			\\
			\left\langle \nabla\cdot   \bu_{h},q_{h}\right\rangle _{\Omega}  &=0&&\forall q_{h}\in S,\label{eq: c}
			\\
			\left\langle \bj_{h},\boldsymbol{e}_{h}\right\rangle _{\Omega} - \left\langle \boldsymbol{B}_{h},\nabla\times\boldsymbol{e}_{h}\right\rangle _{\Omega} &=0&&\forall \boldsymbol{e}_{h}\in \boldsymbol{C},\label{eq: d}
			\\
			\left\langle \partial_{t}\bB_{h},\boldsymbol{b}_{h}\right\rangle _{\Omega}  + \left\langle \nabla\times \bE_{h},\boldsymbol{b}_{h}\right\rangle _{\Omega}  &= 0 &&\forall \boldsymbol{b}_{h}\in \boldsymbol{D},\label{eq: e}
			\\
			\Rm^{-1}\left\langle \bj_{h},\boldsymbol{J}_{h}\right\rangle _{\Omega} -  \left\langle \bE_{h},\boldsymbol{J}_{h}\right\rangle _{\Omega}  - \mathcal{A}\left(  \bu_{h},{\bH_{h}},\boldsymbol{J}_{h}\right)   + \mathsf{h}\mathcal{A}\left(  \bj_{h},{\bH_{h}} ,\boldsymbol{J}_{h}\right)  &= 0 &&\forall \boldsymbol{J}_{h}\in \boldsymbol{C},\label{eq: f}
			\\
			\left\langle \partial _{t}\bH_{h},\boldsymbol{g}_{h} \right\rangle _{\Omega} +  \Rm^{-1}\left\langle \nabla\times \bH_{h},\nabla\times\boldsymbol{g}_{h}\right\rangle _{\Omega}  -\mathcal{A}\left(  {\bu_{h}} ,\bH_{h},\nabla\times\boldsymbol{g}_{h}\right) \quad \label{eq: g}
			\\
			+ \mathsf{h} \mathcal{A}\left( \nabla\times\bH_{h},{\bB_{h}},\nabla\times\boldsymbol{g}_{h}\right) &= 0&& \forall \boldsymbol{g}_{h}\in\boldsymbol{C}_{0},\nonumber
		\end{align}
	\end{subequations}
	subject to initial conditions $\bu_{h}^0$, $\bB_{h}^0$, $\bH_{h}^0$. Homogeneous boundary conditions $P = 0$, $ \boldsymbol{u}\times\boldsymbol{n} = \boldsymbol{0}$, and $\bB\times\bn = \boldsymbol{0}$ are included as natural boundary conditions through boundary integral terms that have vanished in \eqref{eq: a}, \eqref{eq: b}, and \eqref{eq: d}, respectively. For $\bH_{h}$, the boundary condition $\bH\times\bn = \boldsymbol{0}$ is an essential boundary condition and is enforced through space $\boldsymbol{C}_{0}$. For an extension of the form \eqref{Eq: total system} to general boundary conditions, see \ref{App: general BC}. The analysis is sensitive to boundary conditions, and, in this work, we restrict our analysis to \eqref{Eq: total system}, i.e., the formulation for homogeneous boundary conditions \eqref{eq: bc}, only. 
	
	Different to existing dual-field methods employing two weak forms that are dual to each other, see for example \cite[(13a) - (13c), (13d) - (13f)]{ZHANG2022110868} and \cite[(2.4), (2.5)]{MAO2025114130}, the present work uses only one weak form \eqref{Eq: total system} that already involves a pair of dual variables, i.e. $\bB_{h}$ and $\bH_{h}$, and thus is more concise.
	
	\subsection{Incompressibility, magnetic Gauss's law, and conservation of current density}\label{eq: semi discrete properties}
	In \eqref{Eq: total system}, since we have selected $\bu_{h}\in \boldsymbol{D}$, the Hilbert complex \eqref{Eq: discrete complex} ensures that the divergence operator maps $\bu_{h}$ to $S$ surjectively; \eqref{eq: c} guarantees pointwise divergence-free condition for $\bu_{h}$, 
	\begin{equation}\label{eq: div u}
		\nabla\cdot\bu_{h} = 0.
	\end{equation}
	Namely, the incompressibility is preserved by the formulation. Meanwhile, because we are working with a contractible domain $\Omega$ and the image of $\nabla \times$ on $\boldsymbol{C}$ coincides with the kernel of $\nabla \cdot$ on $\boldsymbol{D}$, \eqref{eq: e} ensures that the time derivative of $\bB_{h}$ must be divergence-free, which is saying that, if the initial condition $\bB^{0}_{h}$ is divergence-free, we have magnetic Gauss law pointwise satisfied, i.e.
	\[
	\nabla\cdot\bB_{h} = 0.
	\]
	Due to the same reason, for $\bH_{h}\in\boldsymbol{C}$, if we define $\boldsymbol{\mathsf{j}}_{h}=\nabla\times\bH_{h}$ to be the dual-field of charge density, it satisfies conservation of current density,
	\begin{equation}\label{semi-discrete conservation of current density}
		\nabla\cdot\boldsymbol{\mathsf{j}}_{h}=0,
	\end{equation}
	
	Note that magnetic Gauss's law is only weakly satisfied by $\bH_{h}\in\boldsymbol{C}$. A proof can be found in \cite{decoupledMHD} where the same magnetic induction equation as \eqref{eq: eB} is used. Meanwhile, for $\bj_{h}\in\boldsymbol{C}$ and $\bB_{h}\in \boldsymbol{D}$, Ampere’s law \eqref{eq: mixed Hall MHD HB d} can only be satisfied in a weak sense through integration by parts \eqref{eq: d}. As a result, conservation of current density is only weakly satisfied by $\bj_{h}$. It is a common phenomenon for dual-fied methods that a certain structure is preserved only by one of the dual-fields \cite{ZHANG2022110868, MAO2025114130}.
	
	\subsection{The semi-discrete energy law}
	The  semi-discrete (total) energy is defined as
	\[
	\mathcal{E}_{h} := \mathcal{K}_{h} + \mathcal{M}_{h}
	\]
	where
	\[
	\mathcal{K}_{h}:=\dfrac{1}{2}\left\langle\bu_{h},\bu_{h}\right\rangle_{\Omega},\quad  \mathcal{M}_{h}:=\dfrac{\mathsf{c}}{2}\left\langle\bB_{h},\bB_{h}\right\rangle_{\Omega}
	\]
	are the kinetic part of the fluid and the magnetic part of the electromagnetic field, respectively. 
	
	To analyze the time derivative of the semi-discrete energy for the formulation \eqref{Eq: total system}, we start with \eqref{eq: a}, which holds $\forall \bv_{h}\in \boldsymbol{D}$ and thus must hold for $\bu_{h}\in \boldsymbol{D}$, namely,
	\begin{equation}\label{eq: E0}
		\left\langle \partial _{t}\bu_{h},\bu_{h}\right\rangle _{\Omega} +\mathcal{A}\left( \bw_{h}, \bu_{h},\bu_{h}\right)    + \Rn^{-1}\left\langle \nabla\times \bw_{h} ,\bu_{h}\right\rangle _{\Omega} - \mathsf{c}\mathcal{A}\left(  \bj_{h} ,{\bH_{h}},\bu_{h}\right)  - \left\langle P_{h},\nabla\cdot\bu_{h}\right\rangle _{\Omega}  = \left\langle \bF,\bu_{h}\right\rangle _{\Omega}.
	\end{equation}
	The second term and the fourth term in \eqref{eq: E0} then vanish because of the skew-symmetry of the trilinear operator $\mathcal{A}$ and the incompressibility \eqref{eq: div u}, respectively. As a result, we have
	\begin{equation}\label{eq: partial E f 0}
		\left\langle \partial _{t}\bu_{h},\bu_{h}\right\rangle _{\Omega}    +\Rn^{-1}\left\langle \nabla\times \bw_{h} ,\bu_{h}\right\rangle _{\Omega} - \mathsf{c}\mathcal{A}\left(  \bj_{h} ,{\bH_{h}},\bu_{h}\right) = \left\langle \bF,\bu_{h}\right\rangle _{\Omega}.
	\end{equation}
	We also know from \eqref{eq: b} that, if we select $\boldsymbol{w}_{h}$ to be $\bw_{h}\in\boldsymbol{C}$,
	\[
	\left\langle\bu_{h}, \nabla\times\boldsymbol{\omega}_{h}\right\rangle _{\Omega}=\left\langle \bw_{h},\boldsymbol{\omega}_{h} \right\rangle _{\Omega}.
	\] 
	Inserting this equation into \eqref{eq: partial E f 0} finally leads to an expression of the time derivative of the semi-discrete fluid kinetic energy,
	\begin{equation}\label{eq: partial E f}
		\partial_{t}\mathcal{K}_{h}=\left\langle \partial _{t}\bu_{h},\bu_{h}\right\rangle _{\Omega}    = -\Rn^{-1}\left\langle \bw_{h},\boldsymbol{\omega}_{h} \right\rangle _{\Omega} + \mathsf{c}\mathcal{A}\left(  \bj_{h} ,{\bH_{h}},\bu_{h}\right) + \left\langle \bF,\bu_{h}\right\rangle _{\Omega}.
	\end{equation}
	
	Next, since \eqref{eq: d} holds $\forall \boldsymbol{e}_{h}\in \boldsymbol{C}$, it must hold for $\bE_{h}\in \boldsymbol{C}$, i.e.,
	\[
	\left\langle \bj_{h},\boldsymbol{E}_{h}\right\rangle _{\Omega} - \left\langle \boldsymbol{B}_{h},\nabla\times\boldsymbol{E}_{h}\right\rangle _{\Omega} = 0.
	\]
	Meanwhile, since \eqref{eq: e} holds $\forall \boldsymbol{b}_{h}\in\boldsymbol{D}$, it must holds for $\bB_{h}\in\boldsymbol{D}$, i.e.,
	\[
	\left\langle \partial_{t}\bB_{h},\boldsymbol{B}_{h}\right\rangle _{\Omega}  + \left\langle \nabla\times \bE_{h},\boldsymbol{B}_{h}\right\rangle _{\Omega}  = 0.
	\]
	These two equations together imply
	\begin{equation}\label{eq: Em1}
		\left\langle \partial_{t}\bB_{h},\boldsymbol{B}_{h}\right\rangle _{\Omega}  + \left\langle \bj_{h},\boldsymbol{E}_{h}\right\rangle _{\Omega} = 0.
	\end{equation}
	We now select $\boldsymbol{J}_{h}$ in \eqref{eq: f} to be $\bj_{h}\in\boldsymbol{C}$ and obtain 
	\begin{equation}\label{eq: Em2}
		\Rm^{-1}\left\langle \bj_{h},\boldsymbol{j}_{h}\right\rangle _{\Omega} -  \left\langle \bE_{h},\boldsymbol{j}_{h}\right\rangle _{\Omega}  - \mathcal{A}\left(  \bu_{h},{\bH_{h}},\boldsymbol{j}_{h}\right)   + \mathsf{h}\mathcal{A}\left(  \bj_{h},{\bH_{h}} ,\boldsymbol{j}_{h}\right)  = 0.
	\end{equation}
	Note that $\mathcal{A}\left(  \bj_{h},{\bH_{h}} ,\boldsymbol{j}_{h}\right) =0$ because of the skew-symmetry of the trilinear operator $\mathcal{A}$; Hall effect has no influence to the mangetic energy. Using \eqref{eq: Em1} and \eqref{eq: Em2}, we can obtain an expression for the time derivative of the semi-discrete magnetic energy,
	\begin{equation}\label{eq: Em}
		\partial_{t}\mathcal{M}_{h}=\mathsf{c}\left\langle \partial _{t}\bB_{h},\bB_{h}\right\rangle _{\Omega} = 
		-\mathsf{c}	\Rm^{-1}\left\langle \bj_{h},\boldsymbol{j}_{h}\right\rangle _{\Omega} 
		+ \mathsf{c}\mathcal{A}\left(  \bu_{h},{\bH_{h}},\boldsymbol{j}_{h}\right)  .
	\end{equation}
	With \eqref{eq: partial E f} and \eqref{eq: Em}, we can finally find the time derivative of the semi-discrete energy,
	\begin{equation}\label{eq: discrete E law}
		\partial_{t}\mathcal{E}_{h}=\partial_{t}\mathcal{K}_{h}+\partial_{t}\mathcal{M}_{h}=
		-\Rn^{-1}\left\langle \bw_{h} ,\bw_{h}\right\rangle _{\Omega} 
		-\mathsf{c}	\Rm^{-1}\left\langle \bj_{h},\boldsymbol{j}_{h}\right\rangle _{\Omega}  
		+ \left\langle \bF,\bu_{h}\right\rangle _{\Omega},
	\end{equation}
	which is consistent with the energy law at the continuous level, \eqref{eq: energy law}. Moreover, at the ideal limit, $\Rn,\ \Rm\to\infty$ and $\boldsymbol{f}=\boldsymbol{0}$, \eqref{eq: discrete E law} shows the proposed formulation preserves conservation of energy,
	\[
	\partial_{t}\mathcal{E}_{h} = 0.
	\]
	
	Note that in \eqref{eq: discrete E law} the two trilinear terms, i.e. $\mathsf{c}\mathcal{A}\left(  \bj_{h} ,{\bH_{h}},\bu_{h}\right)$ and $ \mathsf{c}\mathcal{A}\left(  \bu_{h},{\bH_{h}},\boldsymbol{j}_{h}\right) $,  have canceled out. This cancellation implies the exchange between fluid kinetic energy $\mathcal{K}_{h}$ and magnetic energy $\mathcal{M}_{h}$ is exactly captured, which is a key for the derivation of the semi-discrete energy law \eqref{eq: discrete E law}. In more detail, this cancellation (as well as the elimination of the influence on the magnetic energy due to the Hall effect in \eqref{eq: Em2}) is because of the skew-symmetry of the trilinear operator with respect to the first and the third entries. In other words, having $\bB_{h}$ or its dual variable $\bH_{h}$ as the second entry in these terms does not matter; it will result in the same semi-discrete energy law \eqref{eq: discrete E law}. This, in fact, is a major factor that has encouraged us to introduce dual variables for the magnetic field. By doing so, without a violation of the energy law, we additionally obtain conservation of current density \eqref{semi-discrete conservation of current density} and, more importantly, can design a temporal scheme, as will be explained in the next section, that linearizes the fully discrete formulation while maintaining the conservation properties.
	
	We can define a second semi-discrete magnetic energy using $\bH_{h}$ as
	\[
	\widetilde{\mathcal{M}}_{h}:=\dfrac{\mathsf{c}}{2}\left\langle\bH_{h},\bH_{h}\right\rangle_{\Omega}.
	\]
	By selecting $\boldsymbol{g}_{h}$ in \eqref{eq: g} to be $\bH_{h}\in\boldsymbol{C}_{0}$, we have
	\[
	\left\langle \partial _{t}\bH_{h},\bH_{h} \right\rangle _{\Omega} +  \Rm^{-1}\left\langle \nabla\times \bH_{h},\nabla\times\bH_{h}\right\rangle _{\Omega}  -\mathcal{A}\left(  {\bu_{h}} ,\bH_{h},\nabla\times\bH_{h}\right) 
	+ \mathsf{h} \mathcal{A}\left( \nabla\times\bH_{h},{\bB_{h}},\nabla\times\bH_{h}\right) = 0.
	\]
	We can now once again use the skew-symmetry of the trilinear operator to eliminate the fourth term, which leads to an expression of the time derivative of $\widetilde{\mathcal{M}}_{h}$ written as
	\begin{equation} \label{eq: Em H}
		\partial_{t}\widetilde{\mathcal{M}}_{h}=\mathsf{c}\left\langle \partial _{t}\bH_{h},\bH_{h}\right\rangle _{\Omega} = 
		-\mathsf{c}	\Rm^{-1}\left\langle \boldsymbol{\mathsf{j}}_{h},\boldsymbol{\mathsf{j}}_{h}\right\rangle _{\Omega} 
		+ \mathsf{c}\mathcal{A}\left(  \bu_{h},{\bH_{h}},\boldsymbol{\mathsf{j}}_{h}\right)  ,
	\end{equation} 
	where $\boldsymbol{\mathsf{j}}_{h} = \nabla\times\bH_{h} $. \eqref{eq: Em H} shows that $\widetilde{\mathcal{M}}_{h}$ has a comparable expression for the dissipation rate over time as $\mathcal{M}_{h}$, cf. \eqref{eq: Em}. However, $ \boldsymbol{\mathsf{j}}_{h} \in\boldsymbol{D} $ is different from $ \bj_{h}\in\boldsymbol{C}$; they are a pair of dual variables of current density. Therefore, it is not guaranteed that $\mathcal{A}\left(  \bu_{h},{\bH_{h}},\boldsymbol{\mathsf{j}}_{h}\right) $ cancels out $ \mathcal{A}\left(  \bj_{h} ,{\bH_{h}},\bu_{h}\right)$ at the discrete level. As a result, we are not able to derive a semi-discrete energy law as \eqref{eq: discrete E law} for $\mathcal{K}_{h} + \widetilde{\mathcal{M}}_{h}$ from the proposed formulation.
	
	\section{The fully discrete scheme} \label{Sec: discretization}
	Let a set of increasing time instants be
	\begin{equation}\label{eq: ts}
		\left\lbrace t^0, t^{\frac{1}{2}},t^{1},t^{1+\frac{1}{2}},t^{2},t^{2+\frac{1}{2}},\cdots,t^{k-1},t^{t-\frac{1}{2}},t^k,t^{k+\frac{1}{2}},\cdots\right\rbrace,
	\end{equation}
	where $t^{0}=0$, for $k\in\left\lbrace1,2,3,\cdots\right\rbrace$, $t^{k}-t^{k-1}=\varDelta t $ is a positive constant, and $t^{k-\frac{1}{2}}=\dfrac{t^{k-1}+t^{k}}{2}$. We use a superscript to denote the evaluation of a variable at a certain time instant, i.e., for example, 
	\[\bu_{h}^{k}:=\boldsymbol{u}_{h}\left(\boldsymbol{x},t^{k}\right). \]
	
	We now propose a leapfrog-type temporal discretization for \eqref{Eq: total system} and the resulting fully discrete scheme reads as follows: Over the time sequence \eqref{eq: ts}, given $\boldsymbol{f}\in\left[L^2(\Omega)\right]^3$ and $\left( \bu_{h}^{0},\bw_{h}^{0},\bB_{h}^{0},\bj_{h}^{0},\bH_{h}^{\frac{1}{2}}\right)\in \boldsymbol{D}\times\boldsymbol{C}\times\boldsymbol{D}\times\boldsymbol{C}\times\boldsymbol{C}_{0} $, for $k=1,2,3,\cdots$, successively,\vspace{0.2cm}
	
	\noindent(\textbf{\color{red}step 1}) seek $\left(\bu^{k}_{h},\bw^{k}_{h}, P^{k-\frac{1}{2}}_{h}, \bE^{k-\frac{1}{2}}_{h},\bB^{k}_{h}, \bj^{k}_{h}\right)\in \boldsymbol{D}\times \boldsymbol{C}\times S\times \boldsymbol{C}\times \boldsymbol{D}\times \boldsymbol{C}$, such that, $\forall \left(\bv_{h},\boldsymbol{w}_{h}, q_{h}, \boldsymbol{e}_{h},\boldsymbol{b}_{h}, \boldsymbol{J}_{h}\right)\in \boldsymbol{D}\times \boldsymbol{C}\times S\times \boldsymbol{C}\times \boldsymbol{D}\times \boldsymbol{C}$,
	
	\begin{subequations}\label{Eq: system 1}
		\begin{align}
			\left\langle\dfrac{\bu_{h}^{k}-\bu_{h}^{k-1}}{\varDelta t},\bv_{h}\right\rangle _{\Omega} +\mathcal{A}\left( \bw^{k-1}_{h}, \dfrac{\bu^{k-1}_{h}+\bu^{k}_{h}}{2},\bv_{h}\right)  + \Rn^{-1}\left\langle \nabla\times \dfrac{\bw^{k-1}_{h}+\bw^{k}_{h}}{2} ,\bv_{h}\right\rangle _{\Omega}\hspace{-1cm} \label{eq: 1 a}
			\\
			- \mathsf{c}\mathcal{A}\left( \dfrac{\bj^{k-1}_{h}+\bj^{k}_{h}}{2} ,{\bH^{k-\frac{1}{2}}_{h}},\bv_{h}\right)  - \left\langle P^{k-\frac{1}{2}}_{h},\nabla\cdot\bv_{h}\right\rangle _{\Omega}  &= \left\langle \bF^{k-\frac{1}{2}},\bv_{h}\right\rangle _{\Omega},\nonumber\\
			\left\langle \bw^{k}_{h},\boldsymbol{w}_{h} \right\rangle _{\Omega}- \left\langle\bu^{k}_{h}, \nabla\times\boldsymbol{w}_{h}\right\rangle _{\Omega} &=0,\label{eq: 1 b}
			\\
			\left\langle \nabla\cdot   \bu_{h}^{k},q_{h}\right\rangle _{\Omega}  &=0,\label{eq: 1 c}
			\\
			\left\langle \bj^{k}_{h},\boldsymbol{e}_{h}\right\rangle _{\Omega} - \left\langle \boldsymbol{B}^{k}_{h},\nabla\times\boldsymbol{e}_{h}\right\rangle _{\Omega} &=0,\label{eq: 1 d}
			\\
			\left\langle \dfrac{\bB_{h}^{k}-\bB_{h}^{k-1}}{\varDelta t},\boldsymbol{b}_{h}\right\rangle _{\Omega}  + \left\langle \nabla\times \bE^{k-\frac{1}{2}}_{h},\boldsymbol{b}_{h}\right\rangle _{\Omega}  &= 0 ,\label{eq: 1 e}
			\\
			\Rm^{-1}\left\langle \dfrac{\bj^{k-1}_{h}+\bj^{k}_{h}}{2},\boldsymbol{J}_{h}\right\rangle _{\Omega} -  \left\langle \bE^{k-\frac{1}{2}}_{h},\boldsymbol{J}_{h}\right\rangle _{\Omega}     \hspace{5cm} \label{eq: 1 f}\\
			- \mathcal{A}\left(  \dfrac{\bu^{k-1}_{h}+\bu^{k}_{h}}{2},{\bH^{k-\frac{1}{2}}_{h}},\boldsymbol{J}_{h}\right)
			+ \mathsf{h}\mathcal{A}\left(  \dfrac{\bj^{k-1}_{h}+\bj^{k}_{h}}{2},{\bH^{k-\frac{1}{2}}_{h}} ,\boldsymbol{J}_{h}\right)  &= 0 ,\nonumber
		\end{align}
	\end{subequations}
	(\textbf{\color{blue}step 2}) seek $ \bH_{h}^{k+\frac{1}{2}}\in\boldsymbol{C}_{0}$, such that, $\forall \boldsymbol{g}_{h}\in\boldsymbol{C}_{0} $,
	\begin{equation}\label{Eq: system 2}
		\begin{aligned}
			\left\langle \dfrac{\bH_{h}^{k+\frac{1}{2}}-\bH_{h}^{k-\frac{1}{2}}}{\varDelta t},\boldsymbol{g}_{h} \right\rangle _{\Omega} +  \Rm^{-1}\left\langle \nabla\times \dfrac{\bH^{k-\frac{1}{2}}_{h}+\bH^{k+\frac{1}{2}}_{h}}{2},\nabla\times\boldsymbol{g}_{h}\right\rangle _{\Omega}  \hspace{2.5cm}
			\\-\mathcal{A}\left(  {\bu^{k}_{h}} ,\dfrac{\bH^{k-\frac{1}{2}}_{h}+\bH^{k+\frac{1}{2}}_{h}}{2},\nabla\times\boldsymbol{g}_{h}\right) 
			+ \mathsf{h} \mathcal{A}\left(\nabla\times\dfrac{\bH^{k-\frac{1}{2}}_{h}+\bH^{k+\frac{1}{2}}_{h}}{2},{\bB^{k}_{h}},\nabla\times\boldsymbol{g}_{h}\right) &= 0,
		\end{aligned}
	\end{equation}
	untill $t^k=T$, $t^{k+\frac{1}{2}}>T$, or any other critieron is reached.
	
	It is seen that, to start the iterations, we need to know $\left( \bu_{h}^{0},\bw_{h}^{0},\bB_{h}^{0},\bj_{h}^{0},\bH_{h}^{\frac{1}{2}}\right)\in \boldsymbol{D}\times\boldsymbol{C}\times\boldsymbol{D}\times\boldsymbol{C}\times\boldsymbol{C}_{0}$ in advance. Among them,  $\bu_{h}^{0}$ and $ \bB_{h}^{0}$ are given initial conditions, $\bw_{h}^{0}$ and $\bj_{h}^{0} $ can be computed through
	\begin{equation}\label{eq: w0}
		\left\langle \bw^{0}_{h},\boldsymbol{w}_{h} \right\rangle _{\Omega} =\left\langle\bu^{0}_{h}, \nabla\times\boldsymbol{w}_{h}\right\rangle _{\Omega}\quad \forall \boldsymbol{w}_{h}\in\boldsymbol{C},
	\end{equation}
	and
	\begin{equation}\label{eq: j0}
		\left\langle \bj^{0}_{h},\boldsymbol{e}_{h}\right\rangle _{\Omega}  =\left\langle \boldsymbol{B}^{0}_{h},\nabla\times\boldsymbol{e}_{h}\right\rangle _{\Omega}\quad\forall \boldsymbol{e}_{h}\in\boldsymbol{C},
	\end{equation}
	respectively, and $\bH_{h}^{\frac{1}{2}}$ can be computed by solving
	\begin{equation}\label{Eq: system half}
		\begin{aligned}
			2\left\langle \dfrac{\bH_{h}^{\frac{1}{2}}-\bH_{h}^{0}}{\varDelta t},\boldsymbol{g}_{h} \right\rangle _{\Omega} +  \Rm^{-1}\left\langle \nabla\times \dfrac{\bH^{0}_{h}+\bH^{\frac{1}{2}}_{h}}{2},\nabla\times\boldsymbol{g}_{h}\right\rangle _{\Omega}  \hspace{2.5cm}
			\\-\mathcal{A}\left(  {\bu^{0}_{h}} ,\dfrac{\bH^{0}_{h}+\bH^{\frac{1}{2}}_{h}}{2},\nabla\times\boldsymbol{g}_{h}\right) 
			+ \mathsf{h} \mathcal{A}\left(\nabla\times\dfrac{\bH^{0}_{h}+\bH^{\frac{1}{2}}_{h}}{2},{\bB^{0}_{h}},\nabla\times\boldsymbol{g}_{h}\right) &= 0 \qquad \forall \boldsymbol{g}_{h}\in\boldsymbol{C}_{0},
		\end{aligned}
	\end{equation}
	where $\bH_{h}^{0} $ is also a known initial condition.
	Once $\bw_{h}^{0}$, $\bj_{h}^{0} $ and $\bH_{h}^{\frac{1}{2}}$ are computed, regular iterations, i.e. \eqref{Eq: system 1} and \eqref{Eq: system 2} for $k=1,2,\cdots$, can be performed untill $t^{k+\frac{1}{2}}>T$ or any other critieron is reached; see Fig.~\ref{fig: scheme} for an illustration of the overall time-advancing scheme. So far, we can see a first argument of introducing dual variables of the magnetic field for the proposed scheme; it is to linearize the nonlinear terms that involve the magnetic field while maintaining a second-order accuracy for the temporal discretization of them. In more detail, under this leapfrog-type time-advancing scheme, a step can borrow known information that temporally is at its mid-point from the other step. For example, at the $k$th iteration, ({\color{red}step 1}) \eqref{Eq: system 1} first borrows $\bH_{h}^{k-\frac{1}{2}}$, i.e. the output of ({\color{blue}step 2}) \eqref{Eq: system 2} at the $(k-1)$st iteration, and then ({\color{blue}step 2}) \eqref{Eq: system 2} borrows $\bB^{k}_{h}$ and $\bu_{h}^{k}$, outputs of ({\color{red}step 1}) \eqref{Eq: system 1} at the current iteration. By doing so, the employed second-order implicit-mid-point integrator, instead of, for instance, a first-order explicit Euler integrator, can result in linearized terms in the fully discrete formulations.
	The time-advancing scheme introduced in this section is easily reusable for the formulation of general boundary conditions in \ref{App: general BC}.
	
	\begin{figure}[h!]
		\centering
		\[
		\left\lbrace\begin{matrix}\color{red}
			\bu_{h}^{0}\\\color{red}\bB_{h}^{0}\\\color{Blue}\bH_{h}^{0}
		\end{matrix}\right\rbrace
		\stackrel{\eqref{eq: w0}\eqref{eq: j0}\eqref{Eq: system half}}{\longrightarrow}
		\left\lbrace\begin{matrix}\color{red}
			\bw_{h}^{0}\\\color{red}\bj_{h}^{0}\\{\color{Blue}\bH_{h}^{\frac{1}{2}}}
		\end{matrix}\right\rbrace
		\stackrel{\eqref{Eq: system 1}}{\longrightarrow}
		\left\lbrace{\color{Red}\begin{matrix}
				\bu_{h}^{1}\\\bw_{h}^{1}\\\bB_{h}^{1}\\\bj_{h}^{1}
		\end{matrix}}\right\rbrace
		\stackrel{\eqref{Eq: system 2}}{\longrightarrow}
		{\color{Blue}\bH_{h}^{\frac{3}{2}}}
		\stackrel{\eqref{Eq: system 1}}{\longrightarrow}
		\left\lbrace{\color{Red}\begin{matrix}
				\bu_{h}^{2}\\\bw_{h}^{2}\\\bB_{h}^{2}\\\bj_{h}^{2}
		\end{matrix}}\right\rbrace
		\stackrel{\eqref{Eq: system 2}}{\longrightarrow}
		{\color{Blue}\bH_{h}^{\frac{5}{2}}}
		\longrightarrow
		\cdots
		\]
		\caption{An illustration of the proposed time-advancing scheme. First, using initial conditions, we comptue $ \bw_{h}^{0}$,  $ \bj_{h}^{0}$, $ \bH_{h}^{\frac{1}{2}}$ through \eqref{eq: w0}, \eqref{eq: j0}, \eqref{Eq: system half}, respectively. Next, for $k=1,2,3,\cdots$, successively, we compute ({\color{red}step 1}) \eqref{Eq: system 1} at the $k$th integer time-step from $t^{k-1}$ to $t^{k}$ and ({\color{blue}step 2}) \eqref{Eq: system 2} at the $k$th half-integer time-step from $t^{k-\frac{1}{2}}$ to $t^{k+\frac{1}{2}}$ in a leapfrog-type scheme untill $t^{k+\frac{1}{2}}>T$ or any other critieron is reached. For a certain $k$, ({\color{red}step 1}) and ({\color{blue}step 2}) combined is call the $k$th iteration.} 
		\label{fig: scheme}
	\end{figure}
	
	\subsection{Properties of the fully discrete scheme}
	Although linearizing the problem is beneficial, whether it preserves the desired structures is a more important question. In this section, we study it for the proposed leapfrog-type time-advancing scheme.
	
	For the same reasons as in Section~\ref{eq: semi discrete properties}, except that this time in the context of the fully discrete level, we can easily derive that, if $\nabla\cdot\bB_{h}^{0}=0$,
	\begin{equation}\label{eq discrete conservation}
		\nabla\cdot\bu_{h}^{k} = \nabla\cdot\bB_{h}^{k}=\nabla\cdot\boldsymbol{\mathsf{j}}_{h}^{k-\frac{1}{2}}=\nabla\cdot\boldsymbol{\mathsf{j}}_{h}^{k+\frac{1}{2}}=0,\quad k \in\left\lbrace1,2,3,\cdots\right\rbrace,
	\end{equation}
	where, for example, $\boldsymbol{\mathsf{j}}_{h}^{k-\frac{1}{2}}=\nabla\times\bH_{h}^{k-\frac{1}{2}}$. This implies the proposed scheme preserves the incompressibility of the fluid, magnetic Gauss's law, and conservation of current density at the fully discrete level.
	
	Now we study the behavior of the discrete energy for the fully discrete scheme at the $k$th iteration,  $k\in\left\lbrace1,2,3,\cdots\right\rbrace$. The fully discrete energy at the time instant $t^{l}$, $l\in\left\lbrace0,1,2,\cdots\right\rbrace $, is defined as
	\[
	\mathcal{E}_{h}^{l} := \mathcal{K}_{h}^{l} + \mathcal{M}_{h}^{l}
	\]
	where
	\[
	\mathcal{K}_{h}^{l}:=\dfrac{1}{2}\left\langle\bu^{l}_{h},\bu^{l}_{h}\right\rangle_{\Omega},\quad  \mathcal{M}_{h}^{l}:=\dfrac{\mathsf{c}}{2}\left\langle\bB^{l}_{h},\bB^{l}_{h}\right\rangle_{\Omega}
	\]
	are the kinetic part of the fluid and the magnetic part of the electromagnetic field, respectively.
	
	We first take \eqref{eq: 1 a} that holds $\forall \boldsymbol{v}_{h}\in\boldsymbol{D}$ and, thus, must hold for $\dfrac{\bu_{h}^{k-1}+\bu_{h}^{k}}{2}\in\boldsymbol{D}$,
	\begin{equation}\label{eq: discrete E}
		\begin{aligned}
			\left\langle\dfrac{\bu_{h}^{k}-\bu_{h}^{k-1}}{\varDelta t},\dfrac{\bu_{h}^{k-1}+\bu_{h}^{k}}{2}\right\rangle _{\Omega}  + \Rn^{-1}\left\langle \nabla\times \dfrac{\bw^{k-1}_{h}+\bw^{k}_{h}}{2} ,\dfrac{\bu_{h}^{k-1}+\bu_{h}^{k}}{2}\right\rangle _{\Omega}\hspace{-0.8cm}
			\\
			- \mathsf{c}\mathcal{A}\left( \dfrac{\bj^{k-1}_{h}+\bj^{k}_{h}}{2} ,{\bH^{k-\frac{1}{2}}_{h}},\dfrac{\bu_{h}^{k-1}+\bu_{h}^{k}}{2}\right)  &= \left\langle \bF^{k-\frac{1}{2}},\dfrac{\bu_{h}^{k-1}+\bu_{h}^{k}}{2}\right\rangle _{\Omega},
		\end{aligned}
	\end{equation}
	where we have dropped the zeroed terms $\mathcal{A}\left( \bw^{k-1}_{h}, \dfrac{\bu^{k-1}_{h}+\bu^{k}_{h}}{2},\dfrac{\bu_{h}^{k-1}+\bu_{h}^{k}}{2}\right)$ and $\left\langle P^{k-\frac{1}{2}}_{h},\nabla\cdot\dfrac{\bu_{h}^{k-1}+\bu_{h}^{k}}{2}\right\rangle _{\Omega}$ for the skew-symmetry of the trilinear operator and the incompressibility \eqref{eq discrete conservation}, respectively. Note that analyses before \eqref{eq: discrete E} do not request that the initial condition $\bu_{h}^{0}$ is strictly divergence-free (which commonly is the case). And, if $\nabla\cdot \bu_{h}^{0}\neq 0$, \eqref{eq: discrete E} and relations dependent will not be valid for $k=1$. Similarly, from \eqref{eq: 1 b} and \eqref{eq: w0}, we must have
	\[
	\left\langle \dfrac{\bw_{h}^{k-1}+\bw_{h}^{k}}{2} ,\dfrac{\bw_{h}^{k-1}+\bw_{h}^{k}}{2}\right\rangle _{\Omega} = \left\langle \dfrac{\bu_{h}^{k-1}+\bu_{h}^{k}}{2} ,\nabla\times\dfrac{\bw_{h}^{k-1}+\bw_{h}^{k}}{2}\right\rangle _{\Omega}
	\]
	for the linearity of the inner product. With this equation, \eqref{eq: discrete E} can be rewritten as
	\begin{equation*}
		\begin{aligned}
			\left\langle\dfrac{\bu_{h}^{k}-\bu_{h}^{k-1}}{\varDelta t},\dfrac{\bu_{h}^{k-1}+\bu_{h}^{k}}{2}\right\rangle _{\Omega}  + \Rn^{-1}\left\langle \dfrac{\bw_{h}^{k-1}+\bw_{h}^{k}}{2} ,\dfrac{\bw_{h}^{k-1}+\bw_{h}^{k}}{2}\right\rangle _{\Omega}\hspace{-0.8cm}
			\\
			- \mathsf{c}\mathcal{A}\left( \dfrac{\bj^{k-1}_{h}+\bj^{k}_{h}}{2} ,{\bH^{k-\frac{1}{2}}_{h}},\dfrac{\bu_{h}^{k-1}+\bu_{h}^{k}}{2}\right)  &= \left\langle \bF^{k-\frac{1}{2}},\dfrac{\bu_{h}^{k-1}+\bu_{h}^{k}}{2}\right\rangle _{\Omega},
		\end{aligned}
	\end{equation*}
	which further leads to 
	\begin{equation}\label{eq: discrete energy K}
		\frac{\mathcal{K}^{k}_{h}-\mathcal{K}^{k-1}_{h}}{\varDelta t} = -\Rn^{-1}\left\langle \overline{\bw}^{k-\frac{1}{2}} ,\overline{\bw}^{k-\frac{1}{2}}\right\rangle _{\Omega}- \mathsf{c}\mathcal{A}\left( \overline{\bj}^{k-\frac{1}{2}} ,{\bH^{k-\frac{1}{2}}_{h}},\overline{\bu}^{k-\frac{1}{2}}\right) + \left\langle \bF^{k-\frac{1}{2}},\overline{\bu}^{k-\frac{1}{2}}\right\rangle _{\Omega},
	\end{equation}
	where we have employed the overlined notations to denote the averaged discrete variables using the mid-point rule, i.e., 
	$\overline{\boldsymbol{\alpha}}_{h}^{k-\frac{1}{2}} = \dfrac{\boldsymbol{\alpha}^{k-1}_{h}+\boldsymbol{\alpha}^{k}_{h}}{2}$. Notice the difference between $\overline{\boldsymbol{\alpha}}_{h}^{k-\frac{1}{2}}$ and ${\boldsymbol{\alpha}}_{h}^{k-\frac{1}{2}} := \boldsymbol{\alpha}(\boldsymbol{x},t^{k-\frac{1}{2}})$.
	%
	%\begin{equation}\label{eq: mid-point rule}
	%	\overline{\boldsymbol{a}}^{k-\frac{1}{2}}= \dfrac{\boldsymbol{a}^{k-1}_{h}+\boldsymbol{a}^{k}_{h}}{2},\qquad\qquad
	% 	\overline{\boldsymbol{b}}^{k}= \dfrac{\boldsymbol{b}^{k-\frac{1}{2}}_{h}+\boldsymbol{b}^{k+\frac{1}{2}}_{h}}{2}.
	%\end{equation}
	
	Next, taking \eqref{eq: 1 e},  we can obtain
	\begin{equation*}
		\left\langle \dfrac{\bB_{h}^{k}-\bB_{h}^{k-1}}{\varDelta t},\dfrac{\bB_{h}^{k-1}+\bB_{h}^{k}}{2} \right\rangle _{\Omega}  + \left\langle \nabla\times \bE^{k-\frac{1}{2}}_{h},\dfrac{\bB_{h}^{k-1}+\bB_{h}^{k}}{2}\right\rangle _{\Omega}  = 0 
	\end{equation*}
	as $\dfrac{\bB_{h}^{k-1}+\bB_{h}^{k}}{2}\in\boldsymbol{D}$. Moreover, from \eqref{eq: 1 d} and \eqref{eq: j0}, we can arrive at
	\[
	\left\langle \dfrac{\bj_{h}^{k-1}+\bj_{h}^{k}}{2} ,\bE_{h}^{k-\frac{1}{2}}\right\rangle _{\Omega} = \left\langle \dfrac{\bj_{h}^{k-1}+\bj_{h}^{k}}{2} ,\nabla\times\bE_{h}^{k-\frac{1}{2}}\right\rangle _{\Omega}.
	\]
	These two equations togerther reveal
	\begin{equation}\label{eq discrete energy M 1}
		\left\langle \dfrac{\bB_{h}^{k}-\bB_{h}^{k-1}}{\varDelta t},\dfrac{\bB_{h}^{k-1}+\bB_{h}^{k}}{2} \right\rangle _{\Omega} = - \left\langle \dfrac{\bj_{h}^{k-1}+\bj_{h}^{k}}{2} ,\bE_{h}^{k-\frac{1}{2}}\right\rangle _{\Omega}.
	\end{equation}
	Then we revisit \eqref{eq: 1 f} that holds $\forall \boldsymbol{J}_{h}\in\boldsymbol{C} $ and thus must hold for $\dfrac{\bj_{h}^{k-1}+\bj_{h}^{k}}{2}$, which gives
	\begin{equation}\label{eq discrete energy M 2}
		\begin{aligned}
			\Rm^{-1}\left\langle \dfrac{\bj^{k-1}_{h}+\bj^{k}_{h}}{2},\dfrac{\bj_{h}^{k-1}+\bj_{h}^{k}}{2}\right\rangle _{\Omega} -  \left\langle \bE^{k-\frac{1}{2}}_{h},\dfrac{\bj_{h}^{k-1}+\bj_{h}^{k}}{2}\right\rangle _{\Omega}     \hspace{5cm} \\
			- \mathcal{A}\left(  \dfrac{\bu^{k-1}_{h}+\bu^{k}_{h}}{2},{\bH^{k-\frac{1}{2}}_{h}},\dfrac{\bj_{h}^{k-1}+\bj_{h}^{k}}{2}\right)
			+ \mathsf{h}\mathcal{A}\left(  \dfrac{\bj^{k-1}_{h}+\bj^{k}_{h}}{2},{\bH^{k-\frac{1}{2}}_{h}} ,\dfrac{\bj_{h}^{k-1}+\bj_{h}^{k}}{2}\right)  &= 0
		\end{aligned}
	\end{equation}
	With \eqref{eq discrete energy M 1} and \eqref{eq discrete energy M 2}, we can find a discrete time rate of change for the discrete magnetic energy that reads
	\begin{equation}\label{eq discrete energy M}
		\frac{\mathcal{M}^{k}_{h}-\mathcal{M}^{k-1}_{h}}{\varDelta t} = -\mathsf{c}\Rm^{-1}\left\langle \overline{\bj}_{h}^{k-\frac{1}{2}},\overline{\bj}_{h}^{k-\frac{1}{2}}\right\rangle _{\Omega} - 
		\mathsf{c}\mathcal{A}\left(  \overline{\bu}_{h}^{k-\frac{1}{2}},{\bH^{k-\frac{1}{2}}_{h}},\overline{\bj}_{h}^{k-\frac{1}{2}}\right).
	\end{equation}
	Finally, combining \eqref{eq: discrete energy K} and \eqref{eq discrete energy M}, we find the discrete time rate of change for the fully discrete energy, namely a fully discrete energy law,
	\begin{equation}\label{eq: discrete energy law}
		\frac{\mathcal{E}^{k}_{h}-\mathcal{E}^{k-1}_{h}}{\varDelta t} =
		-\Rn^{-1}\left\langle \overline{\bw}_{h}^{k-\frac{1}{2}} ,\overline{\bw}_{h}^{k-\frac{1}{2}}\right\rangle _{\Omega}
		-\mathsf{c}\Rm^{-1}\left\langle 			\overline{\bj}_{h}^{k-\frac{1}{2}},\overline{\bj}_{h}^{k-\frac{1}{2}}\right\rangle _{\Omega}
		+ \left\langle \bF^{k-\frac{1}{2}},\overline{\bu}^{k-\frac{1}{2}}\right\rangle _{\Omega},
	\end{equation}
	wheres terms $-\mathsf{c}\mathcal{A}\left( \overline{\bj}^{k-\frac{1}{2}} ,{\bH^{k-\frac{1}{2}}_{h}},\overline{\bu}^{k-\frac{1}{2}}\right)$ from \eqref{eq: discrete energy K} and  $-
	\mathsf{c}\mathcal{A}\left(  \overline{\bu}_{h}^{k-\frac{1}{2}},{\bH^{k-\frac{1}{2}}_{h}},\overline{\bj}_{h}^{k-\frac{1}{2}}\right)$ from \eqref{eq discrete energy M} have canceled each other due to the skew-symmetry of the trilinear operator; the exchange between the fluid kinetic energy and magnetic energy is exactly captured in the fully discrete scheme. 
	
	We usually have the freedom to define discrete energy and its dissipation law differently at the discrete level using variables at mixed time instants or in modified formats. The discrete energy law \eqref{eq: discrete energy law}, in the author's opinion, is the most canonical counterpart of the continuous version \eqref{eq: discrete E law} under comparable temporal schemes. 
	And we arrive at \eqref{eq: discrete energy law} by construction; it does not depend on exact numerical integrations for computing the inner product and the trilinear terms.
	
	The fully discrete energy law \eqref{eq: discrete energy law} is in line with the semi-discrete energy law \eqref{eq: discrete E law} and the energy law of the strong form \eqref{eq: energy law}. Namely, the proposed scheme exactly captures the dissipation rate in the dissipative regime and preserves conservation of energy in the ideal limit. 
	
	One may doubt the necessity of introducing the dual variable $\bH_{h}$ and its weak evolution equation (i.e. the weak magnetic induction equation) \eqref{eq: g}; solving the Maxwell's part \eqref{eq: d} - \eqref{eq: f} directly on the half-integer time steps, i.e., as the step 2 of a iteration, see \eqref{Eq: system 2}, can also linearize the systems and leads to the same discrete energy law under a similar leapfrog-type time-advancing scheme. This is true, but it solves for $\bE_{h}$, $\bB_{h}$, and $\bj_{h}$ at the half-integer time steps and thus will cause more computational difficulties than solving the single-variable-valued magnetic induction equation. Not to mention that conservation of current density can be only satisfied by the dual variable, $\boldsymbol{\mathsf{j}}_{h}$, of current density.
	
	\subsection{Algebraic form}
	Suppose
	\[\mathrm{dim}\ \boldsymbol{C} = \mathfrak{C},\quad\mathrm{dim}\ \boldsymbol{C}_{0} = \mathfrak{C}_{0},\quad\mathrm{dim}\ \boldsymbol{D} = \mathfrak{D},\quad \text{and}\quad  \mathrm{dim}\ S = \mathfrak{S}.\]
	Let sets 
	\[\left\lbrace\boldsymbol{\sigma}^{1},\boldsymbol{\sigma}^{2},\cdots,\boldsymbol{\sigma}^{\mathfrak{C}}\right\rbrace, \ \left\lbrace\boldsymbol{\sigma}_{0}^{1},\boldsymbol{\sigma}_{0}^{2},\cdots,\boldsymbol{\sigma}_{0}^{\mathfrak{C}_{0}}\right\rbrace,\ \left\lbrace\boldsymbol{\tau}^{1},\boldsymbol{\tau}^{2},\cdots,\boldsymbol{\tau}^{\mathfrak{D}}\right\rbrace, \ \text{and}\  \left\lbrace{\chi}^{1},{\chi}^{2},\cdots,{\chi}^{\mathfrak{S}}\right\rbrace, \]
	be bases of $\boldsymbol{C}$, $\boldsymbol{C}_{0}$, $\boldsymbol{D}$, and $S$, respectively. For elements in these spaces, their column vectors of expansion coefficients under these bases are denoted with the sign $\vec{\cdot}$, for example,
	\[
	\bw_{h}^{k} = \begin{bmatrix}
		\boldsymbol{\sigma}^{1}&
		\boldsymbol{\sigma}^{2}&
		\cdots&
		\boldsymbol{\sigma}^{\mathfrak{C}}
	\end{bmatrix}\vec{\omega}^{k} \in\boldsymbol{C}.
	\]
	
	Now we can rewrite the fully discrete scheme, \eqref{Eq: system 1} and \eqref{Eq: system 2}, in an algebraic form as follows: Over time sequence \eqref{eq: ts}, given $\boldsymbol{f}\in\left[L^2(\Omega)\right]^3$ and $\left(\vec{u}^{0}, \vec{\omega}^{0},\vec{B}^{0},\vec{j}^{0},\vec{H}^{\frac{1}{2}}\right)\in\mathbb{R}^{\mathfrak{D}}\times\mathbb{R}^\mathfrak{C}\times\mathbb{R}^\mathfrak{D}\times\mathbb{R}^\mathfrak{C}\times\mathbb{R}^{\mathfrak{C}_{0}}$, for $k=1,2,3,\cdots$, successively,\vspace{0.2cm}
	
	\noindent(\textbf{\color{red}step 1}) seek $\left(\vec{u}^{k},\vec{\omega}^{k},\vec{P}^{k-\frac{1}{2}}, \vec{E}^{k-\frac{1}{2}},\vec{B}^{k}, \vec{j}^{k}\right)\in \mathbb{R}^{\mathfrak{D}}\times\mathbb{R}^{\mathfrak{C}}\times\mathbb{R}^{\mathfrak{S}}\times \mathbb{R}^{\mathfrak{C}}\times\mathbb{R}^{\mathfrak{D}}\times\mathbb{R}^{\mathfrak{C}}$, such that, $\forall \left(\vec{v},\vec{\omega}, \vec{q}, \vec{e},\vec{b}, \vec{J}\right)\in \mathbb{R}^{\mathfrak{D}}\times\mathbb{R}^{\mathfrak{C}}\times\mathbb{R}^{\mathfrak{S}}\times \mathbb{R}^{\mathfrak{C}}\times\mathbb{R}^{\mathfrak{D}}\times\mathbb{R}^{\mathfrak{C}}$,
	\begin{equation}\label{Eq: algebraic system 1}
		\begin{aligned}
			&\left[\begin{array}{ccc;{2pt/2pt}ccc}
				\frac{1}{\varDelta t}\mathsf{M}_{\boldsymbol{D}} + \frac{1}{2}\boldsymbol{A}^{k-1}_{\bw} &
				\frac{1}{2\Rn}\mathsf{M}_{\boldsymbol{D}}\mathsf{C} & -\mathsf{D}^{\mathsf{T}}\mathsf{M}_{S}& \boldsymbol{0} & \boldsymbol{0} & -\frac{\mathsf{c}}{2}\mathsf{A}_{\bH}^{k-\frac{1}{2}}
				\\
				-\mathsf{C}^{\mathsf{T}}\mathsf{M}_{\boldsymbol{D}} & \mathsf{M}_{\boldsymbol{C}}&\boldsymbol{0}&\boldsymbol{0}&\boldsymbol{0}&\boldsymbol{0}\\
				\mathsf{M}_{S}\mathsf{D}&\boldsymbol{0}&\boldsymbol{0}&\boldsymbol{0}&\boldsymbol{0}&\boldsymbol{0}\\\hdashline[2pt/2pt]
				\boldsymbol{0} & \boldsymbol{0} & \boldsymbol{0} & \boldsymbol{0} & -\mathsf{C}^{\mathsf{T}}\mathsf{M}_{\boldsymbol{D}} & \mathsf{M}_{\boldsymbol{C}}\\
				\boldsymbol{0} & \boldsymbol{0} & \boldsymbol{0} & \mathsf{M}_{D}\mathsf{C} & \frac{1}{\varDelta t}\mathsf{M}_{\boldsymbol{D}} & \boldsymbol{0}\\
				\frac{1}{2}\left( \mathsf{A}_{\bH}^{k-\frac{1}{2}}\right) ^{\mathsf{T}}&\boldsymbol{0}& \boldsymbol{0} & -\mathsf{M}_{\boldsymbol{C}}&\boldsymbol{0} &\frac{1}{2\Rm}\mathsf{M}_{\boldsymbol{C}} + \frac{\mathsf{h}}{2}\mathbb{A}_{\bH}^{k-\frac{1}{2}}
			\end{array}\right]
			\left[\begin{array}{c}
				\vec{u}^{k}\\
				\vec{w}^{k}\\
				\vec{P}^{k-\frac{1}{2}}\\
				\vec{E}^{k-\frac{1}{2}}\\
				\vec{B}^{k}\\
				\vec{j}^{k}\\
			\end{array}\right]
			\\
			&\qquad\qquad\qquad\qquad\qquad\quad =
			\left[\begin{array}{c}
				\left( \frac{1}{\varDelta t}\mathsf{M}_{\boldsymbol{D}} - \frac{1}{2}\boldsymbol{A}^{k-1}_{\bw}\right) \vec{u}^{k-1}
				-\frac{1}{2\Rn}\mathsf{M}_{\boldsymbol{D}}\mathsf{C}\vec{\omega}^{k-1}
				+\frac{\mathsf{c}}{2}\mathsf{A}_{\bH}^{k-\frac{1}{2}}\vec{j}^{k-1} + \underline{f}^{k-\frac{1}{2}}\\
				\boldsymbol{0}\\
				\boldsymbol{0}\\\hdashline[2pt/2pt]
				\boldsymbol{0}\\
				\frac{1}{\varDelta t}\mathsf{M}_{\boldsymbol{D}}\vec{B}^{k-1}\\
				-\frac{1}{2}\left( \mathsf{A}_{\bH}^{k-\frac{1}{2}}\right) ^{\mathsf{T}}\vec{u}^{k-1} - \left( \frac{1}{2\Rm}\mathsf{M}_{\boldsymbol{C}}+\frac{\mathsf{h}}{2}\mathbb{A}_{\bH}^{k-\frac{1}{2}}\right) \vec{j}^{k-1}
			\end{array}\right],
		\end{aligned}
	\end{equation}
	(\textbf{\color{blue}step 2}) seek $ \vec{H}^{k+\frac{1}{2}}\in\mathbb{R}^{\mathfrak{C}_{0}}$, such that, $\forall \vec{g}\in\mathbb{R}^{\mathfrak{C}_{0}} $,
	\begin{equation}\label{Eq: algebraic system 2}
		\begin{aligned}
			&\left( 
			\dfrac{1}{\varDelta t}\mathsf{M}_{\boldsymbol{C}_{0}}
			+\dfrac{1}{2\Rm}\mathsf{C}_{0}^{\mathsf{T}}\mathsf{M}_{\boldsymbol{D}}\mathsf{C}_{0}
			- \dfrac{1}{2}\mathsf{C}_{0}^{\mathsf{T}}\mathsf{A}_{\bu}^{k}
			+\dfrac{\mathsf{h}}{2}\mathsf{C}^{\mathsf{T}}_{0}\boldsymbol{A}^{k}_{\bB}\mathsf{C}_{0}
			\right) 
			\vec{H}^{k+\frac{1}{2}}
			\\&\hspace{5cm} =
			\left( \dfrac{1}{\varDelta t}\mathsf{M}_{\boldsymbol{C}_{0}} 
			-\dfrac{1}{2\Rm}\mathsf{C}_{0}^{\mathsf{T}}\mathsf{M}_{\boldsymbol{D}}\mathsf{C}_{0}
			+ \dfrac{1}{2}\mathsf{C}_{0}^{\mathsf{T}}\mathsf{A}_{\bu}^{k}
			-\dfrac{\mathsf{h}}{2}\mathsf{C}^{\mathsf{T}}_{0}\boldsymbol{A}^{k}_{\bB}\mathsf{C}_{0}
			\right) \vec{H}^{k-\frac{1}{2}},
		\end{aligned}
	\end{equation}
	%\begin{equation}
	%	\begin{aligned}
		%		\left\langle \dfrac{\bH_{h}^{k+\frac{1}{2}}-\bH_{h}^{k-\frac{1}{2}}}{\varDelta t},\boldsymbol{g}_{h} \right\rangle _{\Omega} +  \Rm^{-1}\left\langle \nabla\times \dfrac{\bH^{k-\frac{1}{2}}_{h}+\bH^{k+\frac{1}{2}}_{h}}{2},\nabla\times\boldsymbol{g}_{h}\right\rangle _{\Omega}  \hspace{2.5cm}
		%		\\-\mathcal{A}\left(  {\bu^{k}_{h}} ,\dfrac{\bH^{k-\frac{1}{2}}_{h}+\bH^{k+\frac{1}{2}}_{h}}{2},\nabla\times\boldsymbol{g}_{h}\right) 
		%		+ \mathsf{h} \mathcal{A}\left(\dfrac{\bH^{k-\frac{1}{2}}_{h}+\bH^{k+\frac{1}{2}}_{h}}{2},{\bB^{k}_{h}},\nabla\times\boldsymbol{g}_{h}\right) &= 0,
		%	\end{aligned}
	%\end{equation}
	untill $t^k=T$, $t^{k+\frac{1}{2}}>T$, or any other critieron is reached. 
	
	The involved notations also include the symmetric mass matrices $\mathsf{M}_{\boldsymbol{C}}$, $\mathsf{M}_{\boldsymbol{C}_0}$, $\mathsf{M}_{\boldsymbol{D}}$ and $\mathsf{M}_{S}$ of $\boldsymbol{C}$, $\boldsymbol{C}_{0}$, $\boldsymbol{D}$ and $S$, respectively. For example, the entry at the $i$th row and the $j$th column of $\mathsf{M}_{\boldsymbol{C}}$ is
	\[
	\left.\mathsf{M}_{\boldsymbol{C}}\right|_{i,j} = \left\langle\boldsymbol{\sigma}^{i}, \boldsymbol{\sigma}^{j}\right\rangle_{\Omega},\quad i,j\in\left\lbrace1,2,\cdots,\mathfrak{C}\right\rbrace.
	\]
	Matrices $\mathsf{C}$, $\mathsf{C}_{0}$, $\mathsf{D}$ are the discrete curl and divergence operators on $\boldsymbol{C}$, $\boldsymbol{C}_{0}$ and $\boldsymbol{D}$, respectively. And $\underline{f}^{k-\frac{1}{2}}$ is a column vector whose $i$th entry is
	\[
	\left.\underline{f}^{k-\frac{1}{2}}\right|_{i} = \left\langle\boldsymbol{f}^{k-\frac{1}{2}}, \boldsymbol{\tau}^{i}\right\rangle_{\Omega},\quad i\in\left\lbrace1,2,\cdots,\mathfrak{D}\right\rbrace.
	\]
	The matrices with respect to the trilinear terms are $\boldsymbol{A}_{\bw}^{k-1}$,  $\mathsf{A}_{\bH}^{k-\frac{1}{2}}$, $\mathbb{A}_{\bH}^{k-\frac{1}{2}}$, $\mathsf{A}_{\bu}^{k}$ and $\boldsymbol{A}_{\bB}^{k}$ whose entries can be computed by
	\begin{subequations}
		\begin{align*}
			&\boldsymbol{A}_{\bw}^{k-1} = \mathcal{A}\left( \bw^{k-1}_{h}, \boldsymbol{\tau}^{i},\boldsymbol{\tau}^{j}\right),&&\hspace{-1cm} i,j\in\left\lbrace1,2,\cdots,\mathfrak{D}\right\rbrace,
			\\
			&\mathsf{A}_{\bH}^{k-\frac{1}{2}}= \mathcal{A}\left( \boldsymbol{\sigma}^{i} ,{\bH^{k-\frac{1}{2}}_{h}},\boldsymbol{\tau}^{j}\right),&&\hspace{-1cm} i\in\left\lbrace1,2,\cdots,\mathfrak{C}\right\rbrace,\ j\in \left\lbrace1,2,\cdots,\mathfrak{D}\right\rbrace
			\\
			&\mathbb{A}_{\bH}^{k-\frac{1}{2}}= \mathcal{A}\left(  \boldsymbol{\sigma}^{i},{\bH^{k-\frac{1}{2}}_{h}} ,\boldsymbol{\sigma}^{j}\right),&&\hspace{-1cm} i,j\in\left\lbrace1,2,\cdots,\mathfrak{C}\right\rbrace,
			\\
			&\mathsf{A}_{\bu}^{k} =\mathcal{A}\left(  {\bu^{k}_{h}} ,\boldsymbol{\sigma}^{i}_{0},\boldsymbol{\tau}^{j}\right) ,&&\hspace{-1cm} i\in\left\lbrace1,2,\cdots,\mathfrak{C}_{0}\right\rbrace,\ j\in \left\lbrace1,2,\cdots,\mathfrak{D}\right\rbrace,
			\\
			&\boldsymbol{A}_{\bB}^{k} = \mathcal{A}\left(\boldsymbol{\sigma}^{i}_{0},{\bB^{k}_{h}},\boldsymbol{\tau}^{j}\right),&&\hspace{-1cm} i\in\left\lbrace1,2,\cdots,\mathfrak{C}_{0}\right\rbrace,\ j\in \left\lbrace1,2,\cdots,\mathfrak{D}\right\rbrace.
		\end{align*}
	\end{subequations}
	
	We can also derive the algebraic forms,
	\begin{subequations}
		\begin{align*}
			\mathsf{M}_{\boldsymbol{C}}\vec{\omega}^{0} &= \mathsf{C}^{\mathsf{T}}\mathsf{M}_{\boldsymbol{D}}\vec{u}^{0},
			\\
			\mathsf{M}_{\boldsymbol{C}}\vec{j}^{0}  &= \mathsf{C}^{\mathsf{T}}\mathsf{M}_{\boldsymbol{D}}\vec{B}^{0},
			\\
			\left( 
			\dfrac{2}{\varDelta t}\mathsf{M}_{\boldsymbol{C}_{0}} 
			+ \dfrac{1}{2\Rm}\mathsf{C}_{0}^{\mathsf{T}}\mathsf{M}_{\boldsymbol{D}}\mathsf{C}_{0}
			- \dfrac{1}{2}\mathsf{C}_{0}^{\mathsf{T}}\mathsf{A}_{\bu}^{0} 
			+ \dfrac{\mathsf{h}}{2}\mathsf{C}^{\mathsf{T}}_{0}\boldsymbol{A}_{\bB}^{0}\mathsf{C}_{0}
			\right) \vec{H}^{\frac{1}{2}} &= 
			\left( 
			\dfrac{2}{\varDelta t}\mathsf{M}_{\boldsymbol{C}_{0}} 
			- \dfrac{1}{2\Rm}\mathsf{C}_{0}^{\mathsf{T}}\mathsf{M}_{\boldsymbol{D}}\mathsf{C}_{0}
			+ \dfrac{1}{2}\mathsf{C}_{0}^{\mathsf{T}}\mathsf{A}_{\bu}^{0} 
			- \dfrac{\mathsf{h}}{2}\mathsf{C}^{\mathsf{T}}_{0}\boldsymbol{A}_{\bB}^{0}\mathsf{C}_{0}
			\right) \vec{H}^{0},
		\end{align*}
	\end{subequations}
	of the systems, \eqref{eq: w0} - \eqref{Eq: system half}, respectively, for $\vec{\omega}^{0}$, $\vec{j}^{0}$ and $\vec{H}^{\frac{1}{2}}$ that should be solved ahead of regular iterations subject to known initial conditions $\vec{u}^{0} $, $\vec{B}^{0}$ and $\vec{H}^{0}$. Recall the over time-advancing scheme in Fig.~\ref{fig: scheme}.
	
	We may take the time to reflect upon the algebraic form \eqref{Eq: algebraic system 1}. If we regard \eqref{Eq: algebraic system 1} as a standard $\boldsymbol{Ax}=\boldsymbol{b}$ form,  in its $\boldsymbol{A}$ matrix we can see a clear coupled structure of the Navier-Stokes part and the Maxwell part; the top-left Navier-Stokes part is coupled with the bottom-right Maxwell part through the trilinear terms with respect to $\bj_{h}$, $\bH_{h}$ and $\bu_{h}$ whose algebraic representations now have their places at the top-right and bottom-left corners of $\boldsymbol{A}$. 
	
	\section{Numerical tests} \label{Sec: numerical}
	In this section, we present three numerical tests: a temporal accuracy test, a spatial accuracy test, and a test for the structure-preserving properties. We use mimetic spectral spaces \cite{kreeft2013mixed,palha2014physics,zhang2021hybrid} as our finite-dimensional spaces. Any other set of finite-dimensional spaces that satisfies the discrete Hilbert complex \eqref{Eq: discrete complex} and the regularity \eqref{Eq:regularity} works. The programme is coded in Python. The source code and all data are available on request. 
	
	\subsection{Temporal accuracy test}\label{Sub: temporal test}
	For a given set of parameters $\Rn$, $\Rm$, $\mathsf{c}$ and $\mathsf{h}$, we consider manufactured expressions,
	\[
	\boldsymbol{u}_{1} = \begin{bmatrix}
		\cos x \sin  y \sin  z \ e^{t}\\
		\sin x   \cos y \sin  z \ e^{t}\\
		-2\sin x    \sin  y \cos  z \ e^{t}
	\end{bmatrix},\ 
	\boldsymbol{E}_{1}^{\prime} = \begin{bmatrix}
		\sin x  \cos  y  \ e^{t}\\
		- \sin  y  \cos  z\  e^{t}\\
		-\cos x \sin  z \  e^{t}
	\end{bmatrix},\ 
	\bB^{0}_{1} =\boldsymbol{0},\ \text{and}\ 
	P_{1} = \sin x \sin  y \sin z \ e^{-t}.
	\]
	From them, we compute $\bB_{1}$, $\bj_{1}$ and $\bE_{1}$ using equations \eqref{eq: mixed Hall MHD HB e}, \eqref{eq: mixed Hall MHD HB d} and \eqref{eq: mixed Hall MHD HB f}, i.e.,
	\[\bB_{1} = -\int_{0}^{T}\nabla\times\bE_{1}^\prime\ \mathrm{d}t,\quad \bj_{1}=\nabla\times\bB_{1},\quad \bE_{1} = \Rm^{-1}\bj_{1}-\bu_{1}\times\bB_{1} + \mathsf{h}\bj_{1}\times\bB_{1},\] 
	and we will find that $\bE_{1}\neq\bE_{1}^{\prime}$. To balance the equations, we introduce an extra magnetic source term $\boldsymbol{m}_{1} := \nabla\times\left( \boldsymbol{E}_{1} - \boldsymbol{E}_{1}^{\prime}\right) $ such that
	\[
	\partial_{t}\boldsymbol{B}_{1} + \nabla\times\boldsymbol{E}_{1} = \boldsymbol{m}_{1}.
	\]
	This is common in tests of numerical methods for (Hall) MHD problems (also see \cite[Section~4.1]{doi:10.1137/23M1553844}) since it is usually very difficult to construct a set of manufactured solutions free of a magnetic source term. Discrete formulations with the magnetic source term can be easily derived accordingly. 
	Finally, the manufactured solutions used in this test are
	\begin{equation}\label{eq: manu solution 1}
		\bu_{1},\quad \bw_{1},\quad P_{1},\quad \bE_{1},\quad \bB_{1}, \quad\bj_{1},
	\end{equation}
	subject to source terms, i.e. the external body force $\boldsymbol{f}_{1}$ and the magnetic source $\boldsymbol{m}_{1}$. Expressions of $\bw_{1}$ and $\bF_{1}$ can be calculated through \eqref{eq: mixed Hall MHD HB b} and \eqref{eq: mixed Hall MHD HB a}, respectively. 
	
	The space-time domain is selected to be $\Omega_{1}\times\left(0, T \right] :=\left[0,2\pi\right]^3\times\left( 0, 1\right] $. One can prove that, on this domain, boundary conditions \eqref{eq: bc} and the divergence-free condition of the initial condition $\bB_{1}^{0}$ are satisfied by the manufactured solutions \eqref{eq: manu solution 1}. A mesh of $K^3$ uniform cubic elements is employed, and $h=\frac{2\pi}{K}$ is the element edge length. We use $N$ to indicate the degree of the finite-dimensional spaces. Errors of solutions at the last iteration, i.e. at the $k$th iteration whose $t^{k} = T$, are measured in different norms.
	
\begin{figure}[h!]
		\centering
		\begin{minipage}[c]{1\textwidth}
			\centering{
				\subfloat{ 
					\begin{minipage}[b]{0.325\textwidth}
						\centering
						\includegraphics[width=1\linewidth]{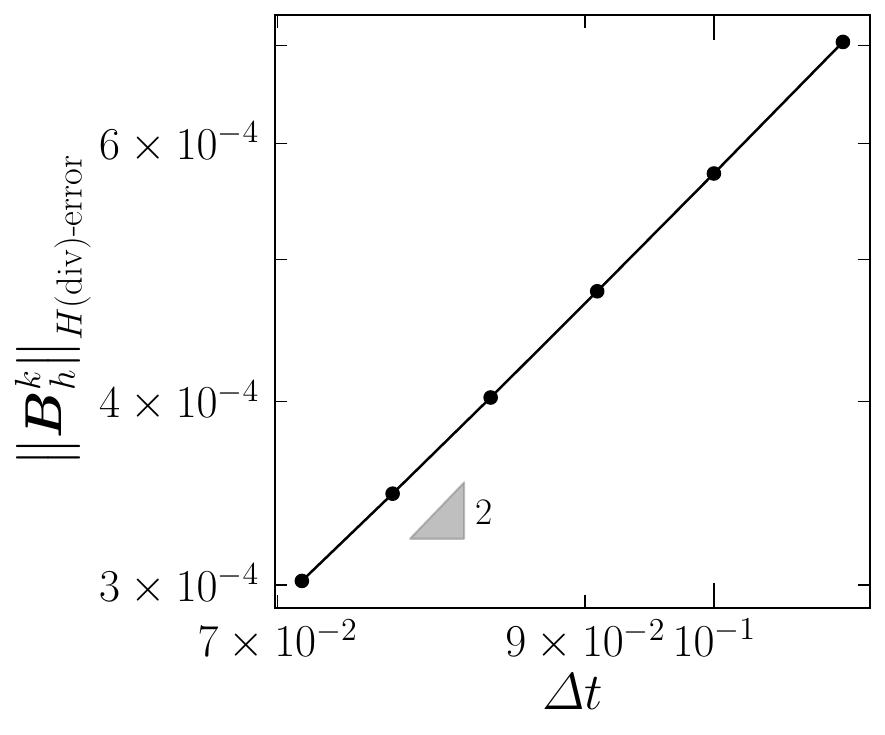}
					\end{minipage}
				}
				\subfloat{
					\begin{minipage}[b]{0.325\textwidth}
						\centering
						\includegraphics[width=1\linewidth]{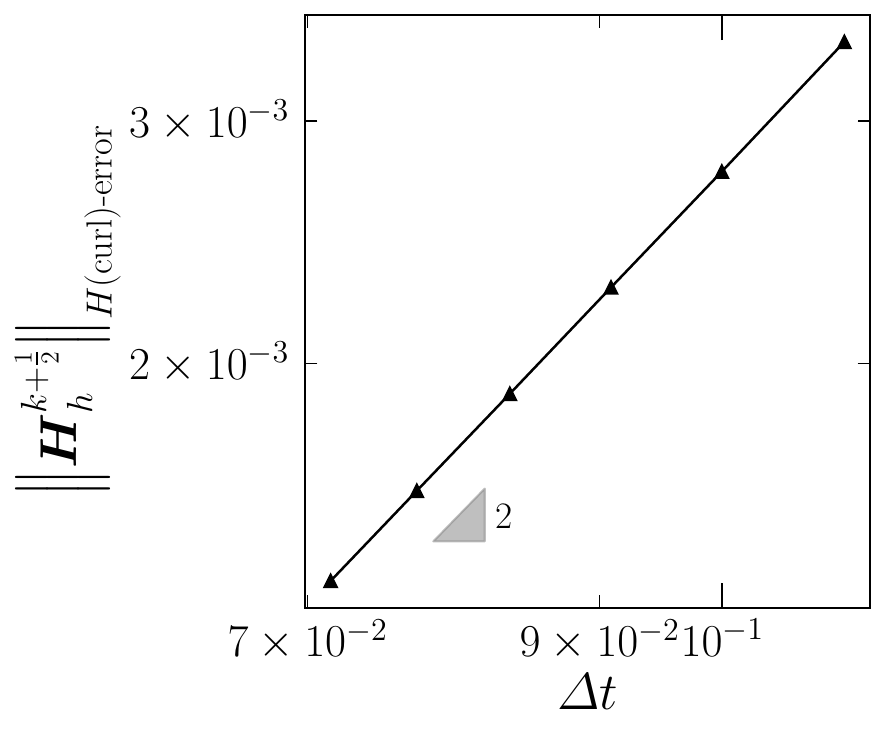}
					\end{minipage}
				}
				\subfloat{
					\begin{minipage}[b]{0.325\textwidth}
						\centering
						\includegraphics[width=1\linewidth]{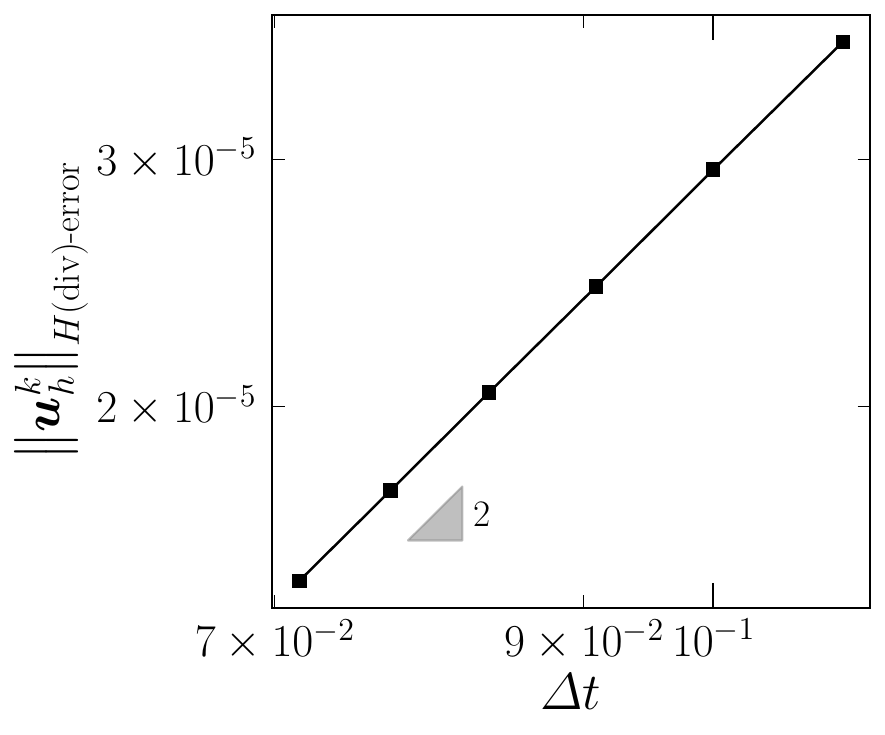}
					\end{minipage}
				}
			}
		\end{minipage}%
		\caption{Results of the temporal accuracy test for $N=3$, $K=6$, $\Rn=\Rm=\mathsf{c}=\mathsf{h}=1$, and $\varDelta t\in\left\lbrace \frac{1}{9}, \frac{1}{10}, \cdots, \frac{1}{14}\right\rbrace$.} 
		\label{fig: results temporal convergence}
	\end{figure}

	In Fig.~\ref{fig: results temporal convergence}, we present some results of the temporal accuracy test for $N=3$, $K=6$, $\Rn=\Rm=\mathsf{c}=\mathsf{h}=1$, and $\varDelta t\in\left\lbrace\frac{1}{9}, \frac{1}{10}, \cdots, \frac{1}{14}\right\rbrace$. These results show that a second-order accuracy is observed for all temporal derivative-related variables, $\boldsymbol{B}_h^k$, $\boldsymbol{H}_h^{k+\frac{1}{2}}$, and $\boldsymbol{u}_h^k$.
	
	\subsection{Spatial accuacy test}
\begin{figure}[h!]
		\centering
		\begin{minipage}[c]{1\textwidth}
			\centering{
				\subfloat{ 
					\begin{minipage}[b]{0.33\textwidth}
						\centering
						\includegraphics[width=1\linewidth]{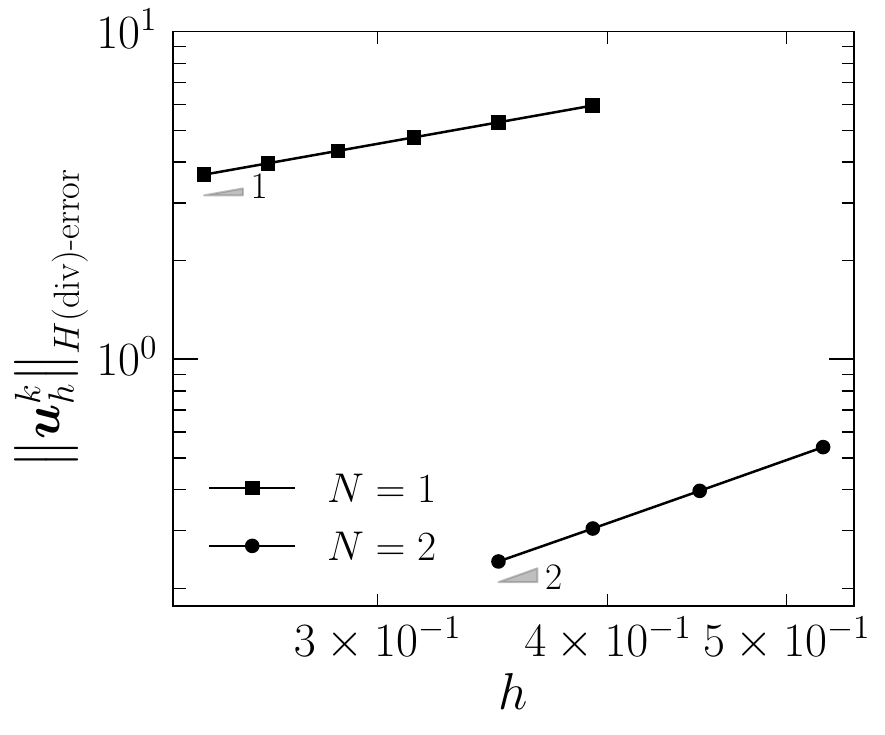}
					\end{minipage}
				}
				\subfloat{
					\begin{minipage}[b]{0.325\textwidth}
						\centering
						\includegraphics[width=1\linewidth]{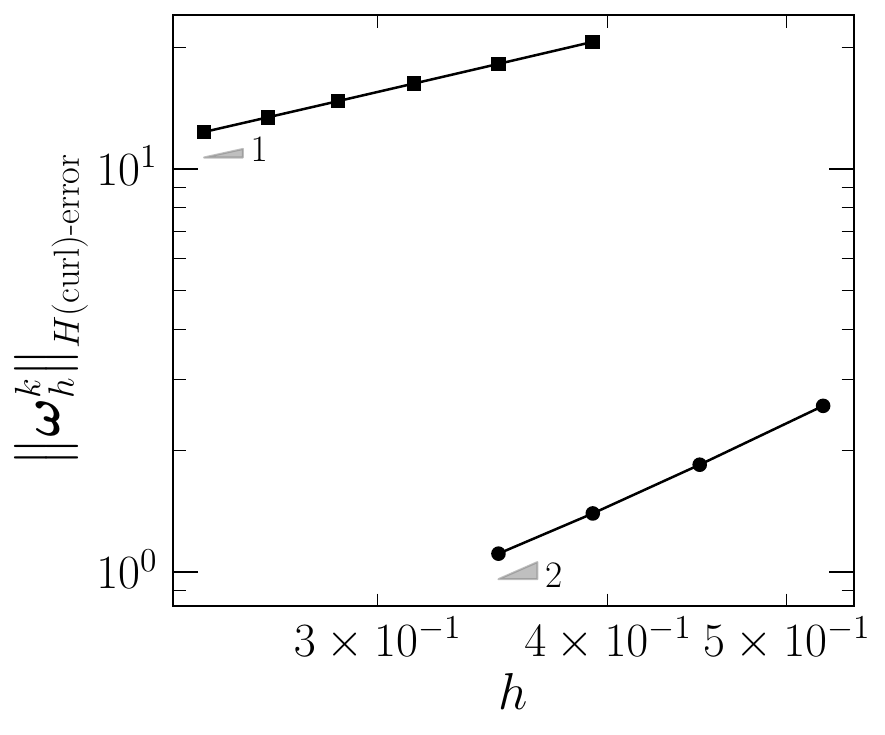}
					\end{minipage}
				}
				\subfloat{
					\begin{minipage}[b]{0.325\textwidth}
						\centering
						\includegraphics[width=1\linewidth]{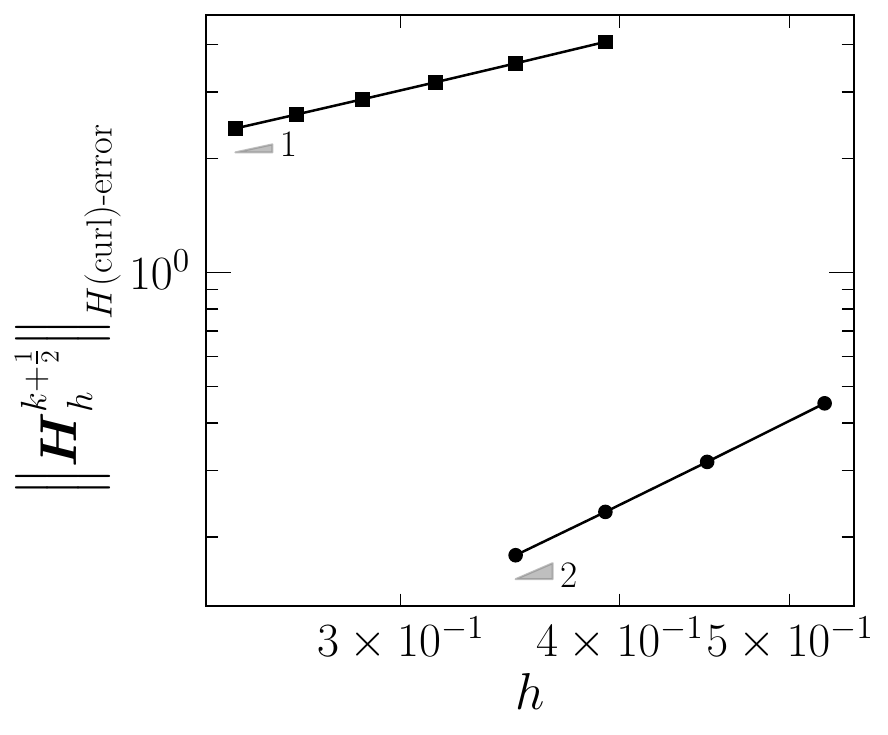}
					\end{minipage}
				}\\
				\subfloat{ 
					\begin{minipage}[b]{0.33\textwidth}
						\centering
						\includegraphics[width=1\linewidth]{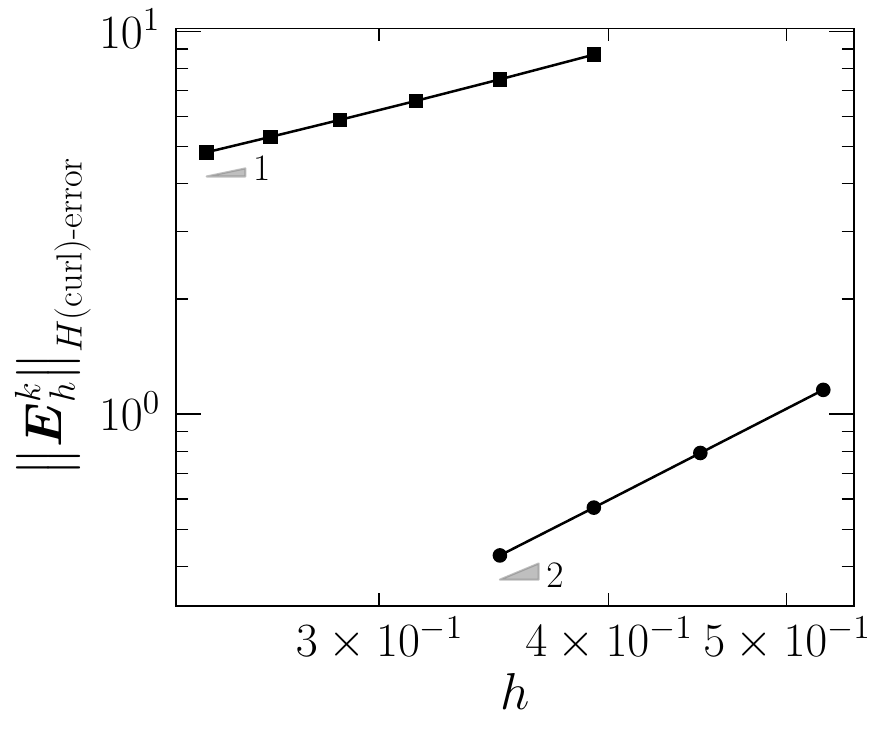}
					\end{minipage}
				}
				\subfloat{
					\begin{minipage}[b]{0.325\textwidth}
						\centering
						\includegraphics[width=1\linewidth]{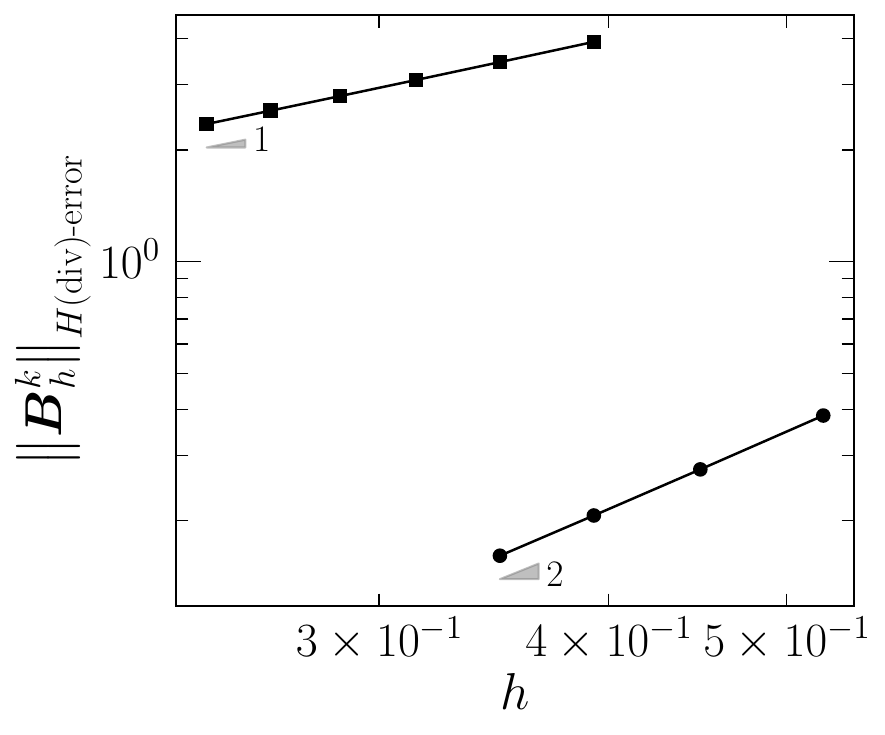}
					\end{minipage}
				}
				\subfloat{
					\begin{minipage}[b]{0.324\textwidth}
						\centering
						\includegraphics[width=1\linewidth]{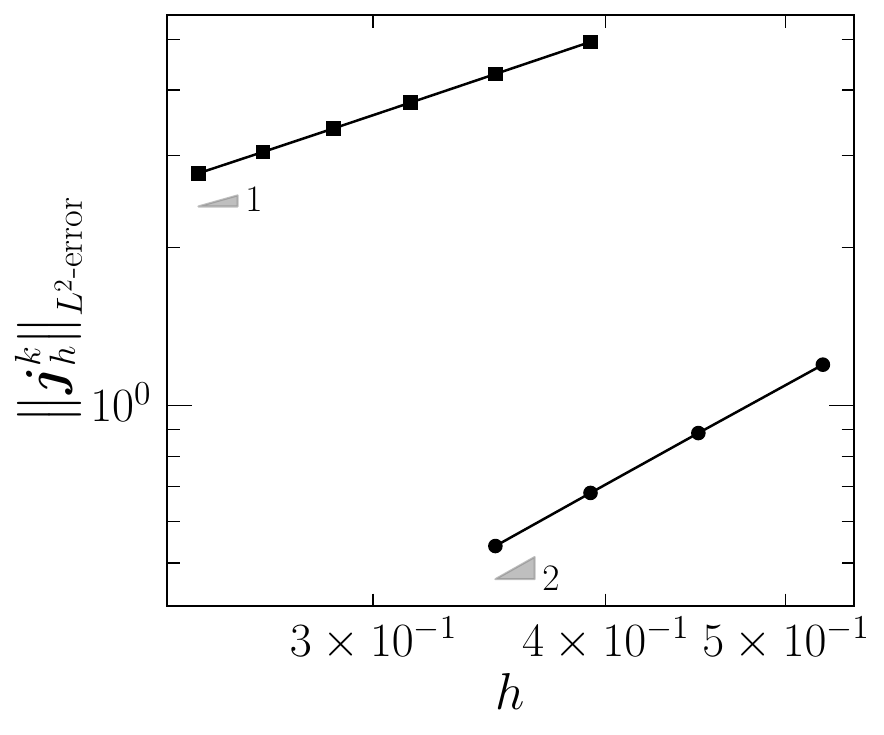}
					\end{minipage}
				}
			}
		\end{minipage}%
		\caption{Results of the spatial accuracy test for $\Rn=\Rm=\mathsf{c}=\mathsf{h}=1$, $\varDelta t=\frac{1}{100}$ and $N\in\left\lbrace1,2\right\rbrace$. For $N=1$, $K\in\left\lbrace16,18,\cdots,26\right\rbrace$. And for $N=2$, $K\in\left\lbrace12,14,16,18\right\rbrace$.} 
		\label{fig: results spatial convergence}
	\end{figure}

	We use the same configuration as in Section~\ref{Sub: temporal test} but with different spatial parameters for this test. Some results for $\Rn=\Rm=\mathsf{c}=\mathsf{h}=1$ and $N\in\left\lbrace1,2\right\rbrace$ under different $h$-refinements are shown in Fig.~\ref{fig: results spatial convergence}. The time interval is selected to be $\varDelta t = \frac{1}{100}$ to avoid the interference of the temporal error. Optimal spatial convergence rates are observed in these results.
	
	\subsection{Structure-preservation test}
\begin{figure}[h!]
		\centering
		\begin{minipage}[c]{1\textwidth}
			\centering{
				\subfloat{ 
					\begin{minipage}[c]{0.5\textwidth}
						\centering
						\includegraphics[width=0.7\linewidth]{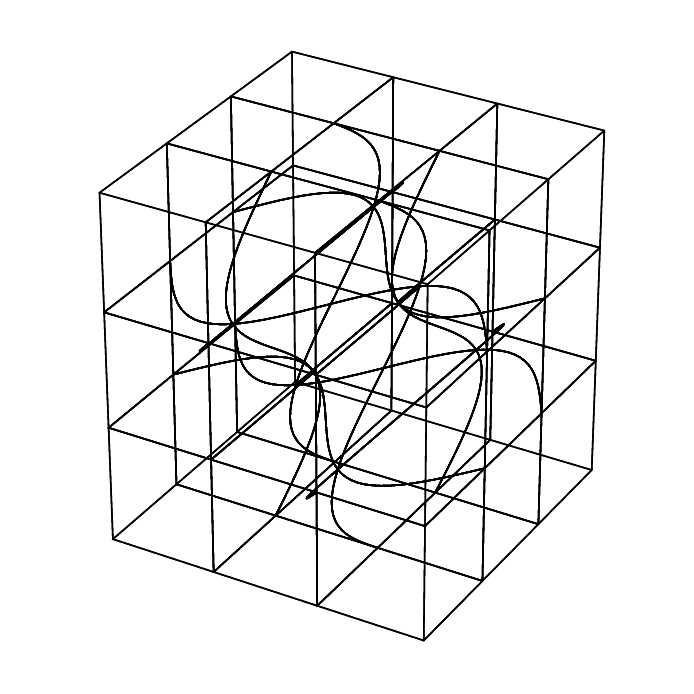}
					\end{minipage}
				}
				\subfloat{
					\begin{minipage}[c]{0.5\textwidth}
						\centering
						\includegraphics[width=0.58\linewidth]{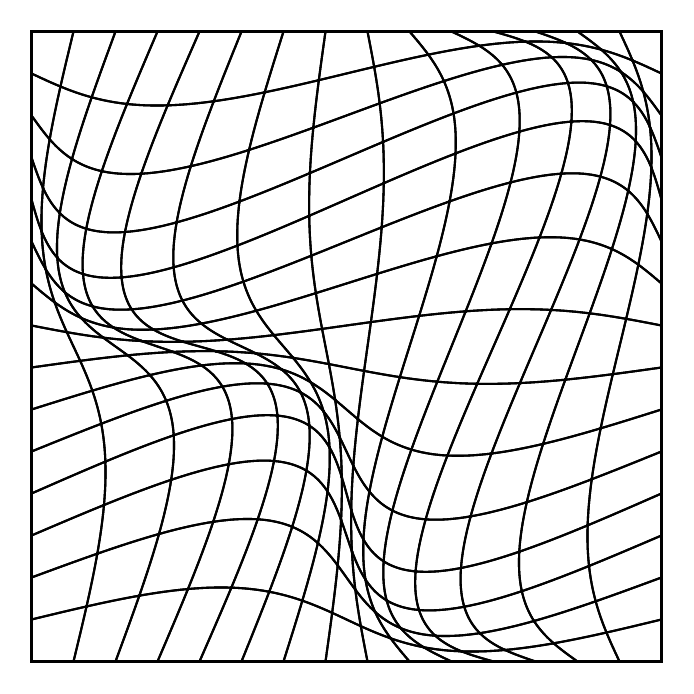}
					\end{minipage}
				}
			}
		\end{minipage}%
		\caption{Left: An illustration of the mesh used in the structure-preservation test for $K=3$ and $d=0.25$. Right: The two-dimensional version of the mesh for $K=15$ that more clearly illustrates the distortion in the mesh.} 
		\label{fig: mesh2d3d}
	\end{figure}

	To test the properties of structure-preservation, we select initial conditions to be
	\begin{equation}\label{eq initial condition 2}
		\bB^{0}_{2} =\boldsymbol{u}^0_{2} = \begin{bmatrix}
			z (z-1)\cos \pi x  \sin  \pi y\\
			z (1-z)\sin \pi x  \cos  \pi y \\
			0
		\end{bmatrix}.
	\end{equation}
	The computational space-time domain is $\Omega_{2}\times\left(0, T\right] :=\left[0,1\right]^3\times\left( 0, 1\right] $. A uniform orthogonal mesh of $K^3$ cubic elements is firstly generated in a reference domain $\Omega_{\mathrm{ref}}:= [0,1]^3$ of a Cartesian coordinate system $(r,s,t)$. This mesh is then transformed into $\Omega_2$ under a nonlinear mapping,
	\begin{equation}\label{crazy mesh}
		\left\lbrace
		\begin{aligned}
			x = r + \dfrac{1}{2}d\sin 2\pi r \sin 2\pi s\sin 2\pi t\\
			y = s + \dfrac{1}{2}d\sin 2\pi r \sin 2\pi s\sin 2\pi t\\
			z = t + \dfrac{1}{2}d\sin 2\pi r \sin 2\pi s\sin 2\pi t
		\end{aligned}\right.,
	\end{equation}
	where $0\leq d\leq 0.25$ indicating the distortion of the mesh. When $d=0$, the mesh is orthogonal. The larger $d$ is, the more distorted the mesh is. An illustration of this mesh for $d=0.25$ can be found in Fig.~\ref{fig: mesh2d3d}, where its two-dimensional version, which more clearly illustrates the distortion of the mesh, is also shown. We use Gauss quadrature for the numerical integration in this paper. The metric behind \eqref{crazy mesh} for $d>0$ cannot be captured by Gauss quadrature; the numerical integration is not exact. With initial conditions \eqref{eq initial condition 2}, on the mesh of $K=9$ and $d=0.25$, the simulation employing the present method is executed under the boundary conditions \eqref{eq: bc} and a zero external body force for parameters $\Rn=\Rm=100$ and $\mathsf{c}=\mathsf{h}=1$ with finite-dimensional spaces of a degree $N=2$. 
	
\begin{figure}[h!]
		\centering
		\begin{minipage}[c]{1\textwidth}
			\centering{
				\subfloat{ 
					\begin{minipage}[c]{0.5\textwidth}
						\centering
						\includegraphics[width=0.9\linewidth]{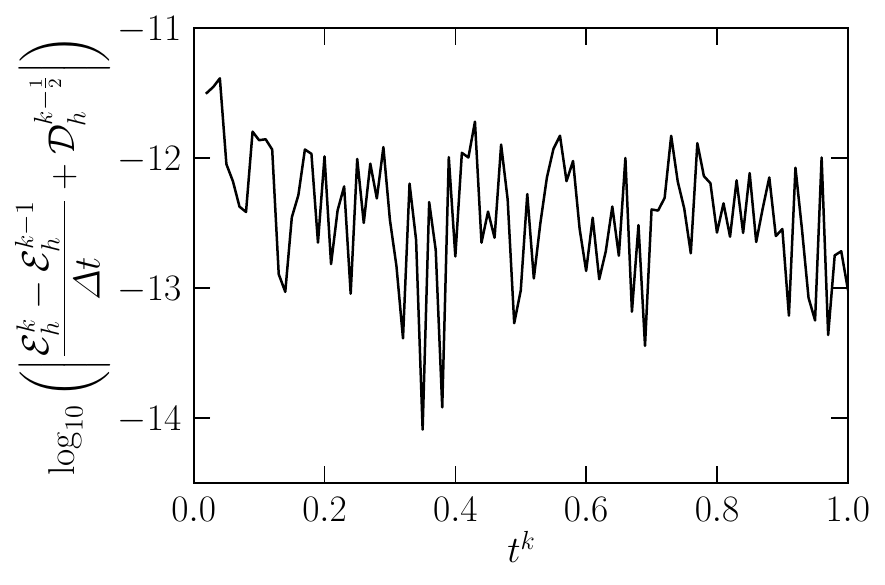}
					\end{minipage}
				}
				\subfloat{
					\begin{minipage}[c]{0.5\textwidth}
						\centering
						\includegraphics[width=0.9\linewidth]{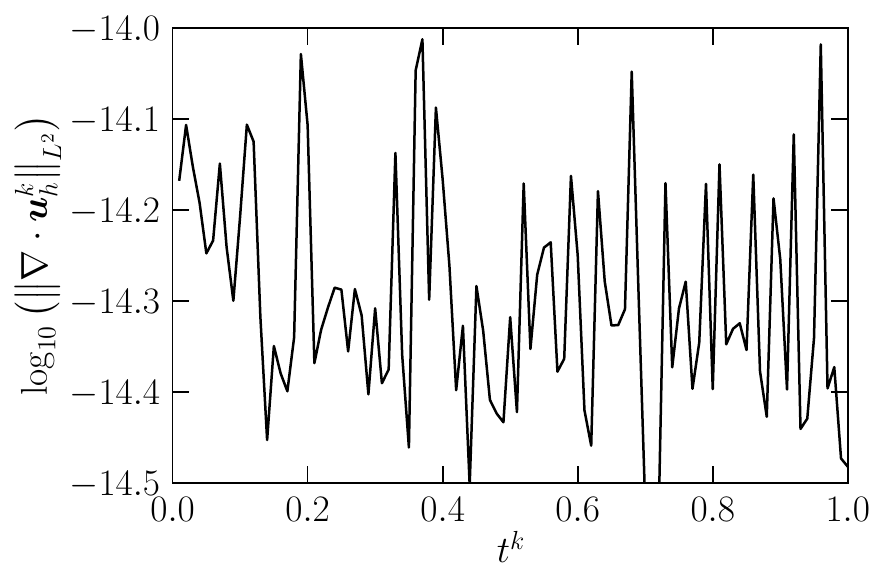}
					\end{minipage}
				}\\
				\subfloat{ 
					\begin{minipage}[c]{0.5\textwidth}
						\centering
						\includegraphics[width=0.9\linewidth]{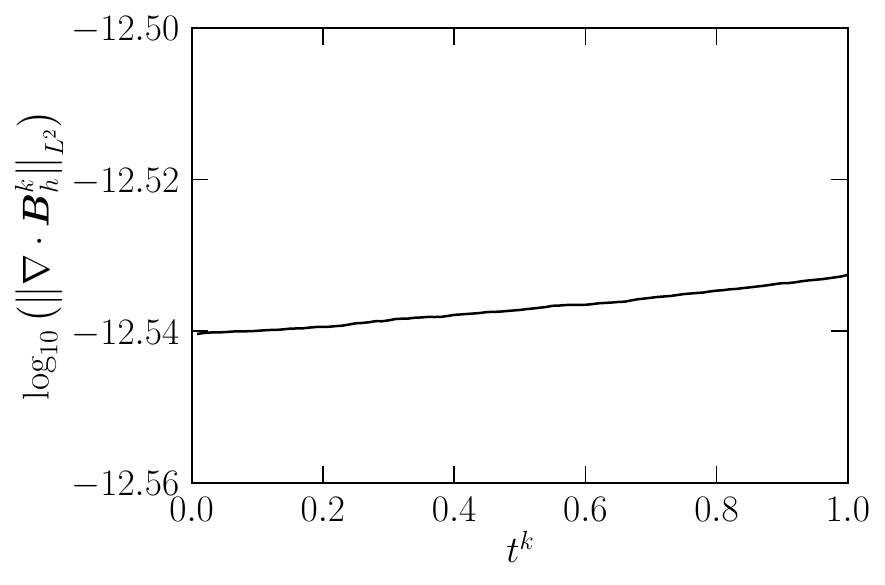}
					\end{minipage}
				}
				\subfloat{
					\begin{minipage}[c]{0.5\textwidth}
						\centering
						\includegraphics[width=0.9\linewidth]{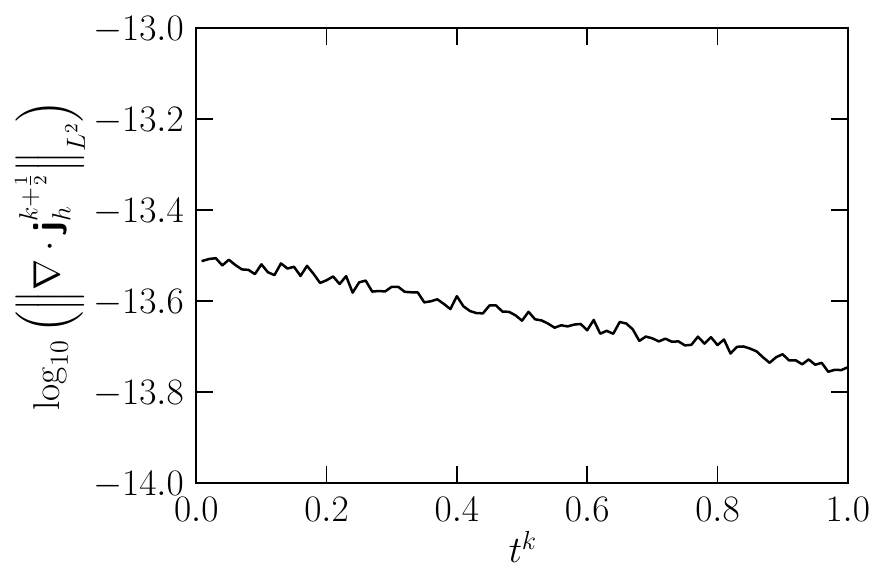}
					\end{minipage}
				}
			}
		\end{minipage}%
		\caption{Some results of the structure-preservation test for $K=9$, $N=2$, $\Rn=\Rm=100$ and $\mathsf{c}=\mathsf{h}=1$. We have introduced $\mathcal{D}_h^{k-\frac{1}{2}}:=\Rn^{-1}\left\langle \overline{\bw}_{h}^{k-\frac{1}{2}} ,\overline{\bw}_{h}^{k-\frac{1}{2}}\right\rangle _{\Omega}
			+\mathsf{c}\Rm^{-1}\left\langle 			\overline{\bj}_{h}^{k-\frac{1}{2}},\overline{\bj}_{h}^{k-\frac{1}{2}}\right\rangle _{\Omega}$ to represent the dissipation rate of energy, c.f. \eqref{eq: discrete energy law}. 
		} 
		\label{fig: conservation results}
	\end{figure}

	In Fig.~\ref{fig: conservation results}, some results of the structure-preservation test are shown. From them, we can see that the discrete energy law \eqref{eq: discrete energy law} and pointwise structures of conservation of mass, magnetic Gauss law, and conservation of current density are all satisfied to the machine precision during the simulation. Due to the fact that no exact numerical integration is used during the simulation, these results support the claim that the present method is structure-preserving by construction. 
	
	\section{Summary and future work} \label{Sec: conclusions}
	This work presents a novel linear structure-preserving mixed finite element discretization of the incompressible Hall MHD model. The scheme satisfies pointwise incompressibility of the fluid, magnetic Gauss's law, and conservation of current density. It also satisfies an energy law by construction such that it captures the energy dissipation rate exactly in the dissipative regime and preserves conservation of energy in the ideal limit. The method is easy to implement in the sense that it only uses the mixed finite element method together with a carefully designed temporal scheme. 
	
	However, the existence, uniqueness, and stability of the solution produced by the proposed method remain unproven. Moreover, although the method is linear, the Navier--Stokes part is not yet decoupled from the Maxwell part of the Hall MHD model. These limitations point to several directions for future research. Further work also includes the extension to conservation of a hybrid helicity.
	
%	However, existence, uniqueness, and stability of the solution for the proposed method are not proven. And, although the method is linear, the Navier-Stokes part is yet not decoupled from the Maxwell part of the Hall MHD model. These limitations suggest the future research of this work. Further work also includes the extension to conservation of a hybrid helicity.
	
	\section*{Acknowledgments}
	This work is supported by the Natural Science Foundation of Guangxi (grant number 2024JJB110005). The author acknowledges dr. Andrea Brugnoli, dr. Artur Palha, dr. Deepesh Toshniwal and dr. Marc Gerritsma for helpful discussions. 
	
	{\small
		\bibliographystyle{elsarticle-num}
		\bibliography{ref}
	}
	
	\appendix
	\section{Differential forms of Hall MHD equations} \label{App: df HMHD}
	The Hall MHD equations, \eqref{eq: mixed Hall MHD}, can be written in differential forms as
	\begin{subequations}\label{eq: df Hall MHD}
		\begin{align}
			\partial _{t}u^{2} +\omega^{1}\wedge\star u^{2}  + \Rn^{-1}\mathrm{d} \omega^{1}  - \mathsf{c}j^{1}\wedge \star B^{2}  + \mathrm{d}^{\ast}P^{3}  &= f^{2},\\
			\omega^{1} &=\mathrm{d}^{\ast}u^{2} ,\\
			\mathrm{d} u^{2}  &=0^{3},\\
			j^{1}  &= \mathrm{d}^{\ast}B^{2}, \label{eq: df Hall MHD j B}\\
			\partial_{t}B^{2} + \mathrm{d}E^{1}  &= 0^{2} ,  \label{eq: df Hall MHD Bt}\\
			\Rm^{-1}j^{1} - \left( E^{1} +\star\left( \star u^{2}\wedge \star B^{2} \right) \right) + \mathsf{h}\star\left(j^{1}\wedge \star B^{2} \right) &= 0^1, \label{eq: df Hall MHD j} 
		\end{align}
	\end{subequations}
	where the variables, i.e. $\bu$, $\bw$, $P$, $\bE$, $\bB$ $\bj$, $\bF$ of \eqref{eq: mixed Hall MHD}, are now expressed as differential forms (or simply forms), $u^2$, $\omega^1$, $P^3$, $E^1$, $B^2$, $j^1$ and $f^2$. The superscript, for example, $k$ of $\alpha^k$, denotes that the form is a $k$-form. Other notations are wedge product $\wedge$, Hodge-$\star$ operator, exterior derivative $\mathrm{d}$, and codifferential $\mathrm{d}^\ast$. We use $0^k$ to express a zero-valued $k$-form.
	
	In three dimensions $(d=3)$, we have, for $k\in\left\lbrace0,1,2,3\right\rbrace$,
	\begin{equation} \label{eq: app star star}
		\star A^k = A^{d-k} ,
		\quad 
		\star\star A^k = A^k,
	\end{equation}
	and, for $l\in\left\lbrace1,2,3\right\rbrace$,
	\begin{equation}\label{eq: app d ast}
		\mathrm{d}^\ast B^l = (-1)^{d\left( l+1\right)  + 1} \star\mathrm{d}\star B^l.
	\end{equation}
	We now introduce the following dual forms,
	\begin{equation}\label{eq: duall variables}
		H^1 = \star B^2,\quad E^2=\star E^1, \quad j^2=\star j^1.
	\end{equation}
	And if we apply a Hodge-$\star$ to \eqref{eq: df Hall MHD Bt}, we obtain
	\[
	\partial_{t}\star B^2 + \star\ \mathrm{d}\star E^2= \star 0^2,
	\]
	i.e., from \eqref{eq: app d ast} and \eqref{eq: duall variables},
	\begin{equation}\label{eq: partial t B2}
		\partial_{t}H^1 + \mathrm{d}^{\ast}\mathsf{E}^{2}  = 0^{1}.
	\end{equation}
	Similarly, we can apply a Hodge-$\star$ to \eqref{eq: df Hall MHD j} and, according to \eqref{eq: app star star} and \eqref{eq: duall variables}, get
	\begin{equation}\label{eq: app E2}
		\Rm^{-1}j^{2} - \left( E^{2} + \left( \star u^{2}\wedge H^1\right) \right) + \mathsf{h}\left(j^{1}\wedge H^1 \right) = 0^2.
	\end{equation}
	We also know from \eqref{eq: df Hall MHD j B}, i.e., $j^1 = \mathrm{d}^\ast B^2$, that
	\[
	j^1 \stackrel{\eqref{eq: app d ast}}{=} \star\ \mathrm{d}\star B^2 \stackrel{\eqref{eq: duall variables}}{=}\star\ \mathrm{d}H^1,
	\]
	and
	\[
	j^2 \stackrel{\eqref{eq: duall variables}}{=} \star j^1 = \star\star\mathrm{d}H^1 \stackrel{\eqref{eq: app star star}}{=} \mathrm{d}H^1.
	\]
	Thus, \eqref{eq: app E2} can be written as 
	\begin{equation}\label{eq: app E2J2}
		E^2 =\Rm^{-1} \mathrm{d}H^1- \star u^2\wedge H^1  + \mathsf{h}\left(\star\mathrm{d}H^1\wedge H^1\right)  .
	\end{equation}
	Combining \eqref{eq: partial t B2} and \eqref{eq: app E2J2} finally leads to an evolution equation for $H^1$,
	\begin{equation}\label{eq: app eH}
		\partial_{t}H^1 + \Rm^{-1}\mathrm{d}^{\ast} \mathrm{d}H^1- \mathrm{d}^{\ast}\left( \star u^2\wedge H^1\right)  + \mathsf{h}\ \mathrm{d}^{\ast}\left(\star\mathrm{d}H^1\wedge H^1\right)  = 0^{1}.
	\end{equation}
	Grouping \eqref{eq: df Hall MHD} and \eqref{eq: app eH} and using $H^1 = \star B^2$, we get a version of \eqref{eq: mixed Hall MHD BH} in  differential forms,
	\begin{subequations}
		\begin{align*}
			\partial _{t}u^{2} +\omega^{1}\wedge\star u^{2}  + \Rn^{-1}\mathrm{d} \omega^{1}  - \mathsf{c}j^{1}\wedge H^1  + \mathrm{d}^{\ast}P^{3}  &= f^{2},\\
			\omega^{1}  &=\mathrm{d}^{\ast}u^{2},\\
			\mathrm{d} u^{2}  &=0^{3},\\
			j^{1} &= \mathrm{d}^{\ast}B^{2} , \\
			\partial_{t}B^{2} + \mathrm{d}E^{1}  &= 0^{2} , \\
			\Rm^{-1}j^{1} - \left( E^{1} +\star\left( \star u^{2}\wedge H^1 \right) \right) + \mathsf{h}\star\left(j^{1}\wedge H^1 \right) &= 0^1 , \\
			\partial_{t}H^1 + \Rm^{-1} \mathrm{d}^{\ast}\mathrm{d}H^1- \mathrm{d}^{\ast}\left( \star u^2\wedge H^1\right)  + \mathsf{h}\ \mathrm{d}^{\ast}\left(\star\mathrm{d}H^1\wedge \star B^2\right)  &= 0^{1}.
		\end{align*}
	\end{subequations}

	\section{The formulation for general boundary conditions} \label{App: general BC}
	Let $\left\lbrace \Gamma_{\hat{P}}, \Gamma_{\hat{u}} \right\rbrace$, $\left\lbrace \Gamma_{\widehat{\bu}}, \Gamma_{\widehat{\boldsymbol{\omega}}} \right\rbrace$, and $\left\lbrace \Gamma_{\widehat{\bB}}, \Gamma_{\widehat{\bE}} \right\rbrace$ be three partitions of the complete boundary $\partial\Omega$, i.e., for example, $\overline{\Gamma_{\hat{P}}} \cup \overline{\Gamma_{\hat{u}}}=\partial\Omega$, $\Gamma_{\hat{P}}\cap \Gamma_{\hat{u}} = \emptyset $. A general setting of boundary conditions for the Hall MHD equations \eqref{eq: mixed Hall MHD BH} is 
	\begin{equation}\label{eq: general bc}
		\left\lbrace
		\begin{aligned}
			P &= \hat{P}\ &&\text{on } \Gamma_{\hat{P}}\\
			\bu\cdot\bn &= \hat{u} \ &&\text{on }\Gamma_{\hat{u}}
		\end{aligned}
		\right.,\qquad
		\left\lbrace
		\begin{aligned}
			\bu\times\bn &= \widehat{\bu} \ &&\text{on } \Gamma_{\widehat{\bu}}\\
			\bw\times\bn &= \widehat{\bw} \ &&\text{on } \Gamma_{\widehat{\bw}}
		\end{aligned}
		\right.,\qquad
		\left\lbrace
		\begin{aligned}
			\bB\times\bn=\bH\times\bn&= \widehat{\bB}\ &&\text{on } \Gamma_{\widehat{\bB}}\\
			\bE\times\bn &= \widehat{\bE} \ &&\text{on }\Gamma_{\widehat{\bE}}
		\end{aligned}
		\right..
	\end{equation}
	When $\Gamma_{\hat{u}} = \Gamma_{\widehat{\bw}} = \Gamma_{\widehat{\bE}} = \emptyset$, $\hat{P} = 0$, and $\widehat{\bu}=\widehat{\bB} = \boldsymbol{0}$, it becomes the particular configuration \eqref{eq: bc}.
	
	Under the general setting of boundary conditions \eqref{eq: general bc}, and let $\mathcal{T}$ be the trace operator, the mixed weak formulation \eqref{Eq: total system} is now given in a version for general boundary conditions as follows.
	Given $\bF_{h}\in \boldsymbol{D}$ and boundary conditions $\hat{P} \in \mathcal{T}H^1(\Omega,\Gamma_{\hat{P}})$, $\hat{u} \in \mathcal{T}H(\Omega;\mathrm{div}, \Gamma_{\hat{P}})$, $\widehat{\bu}\in \mathcal{T}H(\Omega;\mathrm{curl},\Gamma_{\widehat{\bu}})$, $\widehat{\bw}\in \mathcal{T}H(\Omega;\mathrm{curl},\Gamma_{\widehat{\bu}})$, $ \widehat{\bB}\in \mathcal{T}H(\Omega;\mathrm{curl},\Gamma_{\widehat{\bB}})$, and $\widehat{\bE}\in \mathcal{T}H(\Omega;\mathrm{curl},\Gamma_{\widehat{\bE}})$, seek $\left(\bu_{h},\bw_{h}, P_{h}, \bE_{h},\bB_{h}, \bj_{h},\bH_{h}\right)\in \boldsymbol{D}\times \boldsymbol{C}\times S\times \boldsymbol{C}\times \boldsymbol{D}\times \boldsymbol{C}\times\boldsymbol{C}$, such that 
	\begin{subequations}\label{Eq: total system general BC}
		\begin{align*}
			\left\langle \partial _{t}\bu_{h},\bv_{h}\right\rangle _{\Omega} +\mathcal{A}\left(  \bw_{h}, \bu_{h},\bv_{h}\right)   + \Rn^{-1}\left\langle \nabla\times \bw_{h} ,\bv_{h}\right\rangle _{\Omega}\quad \ 
			\\
			- \mathsf{c}\mathcal{A}\left( \bj_{h} ,{\bH_{h}},\bv_{h} \right) - \left\langle P_{h},\nabla\cdot\bv_{h}\right\rangle _{\Omega}  &= \left\langle \bF_{h},\bv_{h}\right\rangle _{\Omega} - \left<\hat{P},\bv_{h}\cdot\bn\right>_{\Gamma_{\hat{P}}},&&\forall \bv_{h}\in \boldsymbol{D},\nonumber\\
			\left\langle \bw_{h},\boldsymbol{w}_{h} \right\rangle _{\Omega}- \left\langle\bu_{h}, \nabla\times\boldsymbol{w}_{h}\right\rangle _{\Omega} &=-\left\langle\widehat{\bu},\boldsymbol{w}_{h}\right\rangle_{\Gamma_{\widehat{\bu}}},&&\forall \boldsymbol{w}_{h}\in \boldsymbol{C},
			\\
			\left\langle \nabla\cdot\bu_{h},q_{h}\right\rangle _{\Omega}  &=0,&&\forall q_{h}\in S,
			\\
			\left\langle \bj_{h},\boldsymbol{e}_{h}\right\rangle _{\Omega} - \left\langle \boldsymbol{B}_{h},\nabla\times\boldsymbol{e}_{h}\right\rangle _{\Omega} &=-\left\langle\widehat{\bB}, \boldsymbol{e}_{h}\right\rangle_{\Gamma_{\widehat{\bB}}},&&\forall \boldsymbol{e}_{h}\in \boldsymbol{C},
			\\
			\left\langle \partial_{t}\bB_{h},\boldsymbol{b}_{h}\right\rangle _{\Omega}  + \left\langle \nabla\times \bE_{h},\boldsymbol{b}_{h}\right\rangle _{\Omega}  &= 0 ,&&\forall \boldsymbol{b}_{h}\in \boldsymbol{D},
			\\
			\Rm^{-1}\left\langle \bj_{h},\boldsymbol{J}_{h}\right\rangle _{\Omega} -  \left\langle \bE_{h},\boldsymbol{J}_{h}\right\rangle _{\Omega}  -\mathcal{A}\left( \bu_{h},{\bH_{h}},\boldsymbol{J}_{h}\right)   + \mathsf{h}\mathcal{A}\left( \bj_{h},{\bH_{h}} ,\boldsymbol{J}_{h}\right)  &= 0 ,&&\forall \boldsymbol{J}_{h}\in \boldsymbol{C},
			\\
			\left\langle \partial _{t}\bH_{h},\boldsymbol{g}_{h} \right\rangle _{\Omega} +  \Rm^{-1}\left\langle \nabla\times \bH_{h},\nabla\times\boldsymbol{g}_{h}\right\rangle _{\Omega}  -\mathcal{A}\left({\bu_{h}} ,\bH_{h},\nabla\times\boldsymbol{g}_{h}\right) \quad 
			\\
			+ \mathsf{h}\mathcal{A}\left(\nabla\times\bH_{h},{\bB_{h}},\nabla\times\boldsymbol{g}_{h}\right) &= \left\langle\widehat{\bE}, \boldsymbol{g}_{h}\right\rangle_{\Gamma_{\widehat{\bE}}},&& \forall \boldsymbol{g}_{h}\in\boldsymbol{C},\nonumber
		\end{align*}
	\end{subequations}
	subject to boundary conditions $ \bu_{h}\cdot\bn = \hat{u}$ on $\Gamma_{\hat{u}}$, $\bw\times\bn = \widehat{\bw}$ on $\Gamma_{\widehat{\bw}}$, $\bE\times\bn=\widehat{\bE}$ on $\Gamma_{\widehat{\bE}}$, and $\bH\times\bn=\widehat{\bB}$ on $\Gamma_{\widehat{\bB}}$, and initial conditions $\bu^0_{h}$, $\bB_{h}^0$, and $\bH^0_{h}$. It is clear that, when $\Gamma_{\hat{u}} = \Gamma_{\widehat{\bw}} = \Gamma_{\widehat{\bE}} = \emptyset$, $\hat{P} = 0$, and $\widehat{\bu}=\widehat{\bB} = \boldsymbol{0}$, it reduces to the formulation \eqref{Eq: total system}. 
	
\end{document}